\documentclass[12pt,a4paper]{amsart}
\usepackage{a4wide}
\usepackage[english]{babel}
\usepackage[T2A]{fontenc}
\usepackage[utf8]{inputenc} 
\usepackage{amsfonts}
\usepackage{amssymb, amsthm, amscd}
\usepackage{amsmath}
\usepackage{mathtools}
\usepackage{needspace}
\usepackage{etoolbox}
\usepackage{lipsum}
\usepackage{comment}
\usepackage{cmap}
\usepackage[pdftex]{graphicx}
\usepackage[unicode]{hyperref}
\usepackage[matrix,arrow,curve]{xy}
\usepackage[usenames,dvipsnames]{xcolor}
\usepackage{colortbl}
\usepackage{textcomp}
\usepackage{cite}
\usepackage{euscript}
\usepackage{epigraph}

\newcommand{\Z}{\mathbb{Z}}

\newcommand{\N}{\mathbb{N}}

\newcommand{\M}{\mathfrak{M}}
\newcommand{\A}{\mathfrak{A}}

\newcommand{\I}{\mathfrak{I}}
\renewcommand{\u}{\mathfrak{u}}
\renewcommand{\v}{\mathfrak{v}}

\newcommand{\ovl}{\overline}

\newcommand{\eps}{\varepsilon}
\newcommand{\ph}{\varphi}

\newcommand{\sub}{\subseteq}

\renewcommand{\ge}{\geqslant}
\renewcommand{\le}{\leqslant}
\newcommand{\sm}{\setminus}
\newcommand{\map}[3]{#1\colon #2\to #3}
\newcommand{\phan}{\vphantom{,}}

\newcommand{\<}{\langle}
\renewcommand{\>}{\rangle}

\newcommand{\emp}{\varnothing}
\newcommand{\SR}{\mathrm{SR}}
\newcommand{\ASR}{\mathrm{ASR}}

\newcommand{\AFSR}{\mathrm{AFSR}}

\newcommand{\SL}{\mathop{\mathrm{SL}}\nolimits}
\newcommand{\GL}{\mathop{\mathrm{GL}}\nolimits}

\newcommand{\D}{\mathop{\mathrm{D}}\nolimits}
\newcommand{\tc}{\text{,}}
\newcommand{\tp}{\text{.}}
\renewcommand{\tilde}{\widetilde}
\renewcommand{\hat}{\widehat}
\newcommand{\vpi}{\varpi}

\DeclareMathOperator{\Max}{Max}
\DeclareMathOperator{\Um}{Um}
\DeclareMathOperator{\Orb}{Orb}
\DeclareMathOperator{\Eq}{Eq}
\DeclareMathOperator{\BSdim}{BSdim}
\DeclareMathOperator{\Jdim}{Jdim}

\DeclareMathOperator{\Ann}{Ann}
\DeclareMathOperator{\res}{res}

\theoremstyle{plain}
\newtheorem{thm}{Theorem}[section]
\newtheorem{lem}[thm]{Lemma}
\newtheorem{prop}[thm]{Proposition}

\theoremstyle{definition}
\newtheorem{defn}[thm]{Definition}

\theoremstyle{remark}
\newtheorem{rem}[thm]{Remark}

\AtBeginEnvironment{thm}{\begin{samepage}}
	\AtEndEnvironment{thm}{\end{samepage}}
\AtBeginEnvironment{lem}{\begin{samepage}}
	\AtEndEnvironment{lem}{\end{samepage}}
\AtBeginEnvironment{st}{\begin{samepage}}
	\AtEndEnvironment{st}{\end{samepage}}
\AtBeginEnvironment{crit}{\begin{samepage}}
	\AtEndEnvironment{crit}{\end{samepage}}
\AtBeginEnvironment{ax}{\begin{samepage}}
	\AtEndEnvironment{ax}{\end{samepage}}
\AtBeginEnvironment{defn}{\begin{samepage}}
	\AtEndEnvironment{defn}{\end{samepage}}
\AtBeginEnvironment{cor}{\begin{samepage}}
	\AtEndEnvironment{cor}{\end{samepage}}
\AtBeginEnvironment{note}{\begin{samepage}}
	\AtEndEnvironment{note}{\end{samepage}}
\AtBeginEnvironment{prop}{\begin{samepage}}
	\AtEndEnvironment{prop}{\end{samepage}}

\makeatletter
\def\@settitle{\begin{center}%
		\baselineskip14\p@\relax
		\bfseries
		\@title
	\end{center}%
}

\def\@evenhead{\hfil\sc Pavel Gvozdevsky\hfil}
\def\@oddhead{\hfil\sc Bounded reduction\hfil}
\makeatother

\title{Bounded reduction for Chevalley groups of types $E_6$ and $E_7$}

\keywords{Chevalley groups, surjective stability of $K_1$, bounded reduction, polynomial rings}
\subjclass[2020]{19B14(primary), 19B10, 20G35 and 20G41(secondary)} 

\author{Pavel Gvozdevsky*}
\date{}
\address{Department of Mathematics, Bar-Ilan University, 5290002 Ramat Gan, ISRAEL}
\thanks{The paper is written as part of the author's post-doctoral fellowship at Bar-Ilan University, Department of Mathematics. Research is supported by Russian Science Foundation grant (project \textnumero 22-21-00257).}
\email{gvozdevskiy96@gmail.com}

\begin{document}

\maketitle

\epigraph{Dedicated to the memory of the brilliant mathematician Irina Suprunenko}

\begin{abstract}
We prove that an element from the Chevalley group of type $E_6$ or $E_7$ over a polynomial ring with coefficients in a small-dimensional ring can be reduced to an element of certain proper subsystem subgroup by a bounded number of elementary root elements. The bound is given explicitly. This result is an effective version of the early stabilisation of the corresponding $K_1$-functor. We also give part of the proof of similar hypothesis for~$E_8$.
\end{abstract}

\section{Introduction}

This paper deals with the bounded reduction in exceptional Chevalley groups over polynomial rings.

Chevalley groups over certain rings have bounded width with respect to the elementary generators. For example this holds for Dedekind domains of arithmetic type, see \cite{CarterKeller}, \cite{CarterKellerZ}, \cite{Morris}, \cite{TavgenChevalley}, \cite{TavgenTwisted},\cite{MorganRapinchukSury},\cite{ErovenkoRapinchuk2001},\cite{ErovenkoRapinchuk2006}, \cite{VavKunPlot}. Results on such {\it bounded generation} are of great value, for example they are connected to the {\it congruence subgroup property}, see \cite{Lubotzky},\cite{PlatonovRapinchuk}; to Margulis--Zimmer conjecture, see \cite{ShalomWillis}; and have applications in studying strong boundedness, see \cite{TrostSympl}, \cite{TrostChevalley}, \cite{TrostSL2}, \cite{TrostChevalleyStrong}, \cite{TrostStabBGandStrong}. However bounded generation occurs very rarely in the sense that classes of rings for which it is known to hold are pretty narrow. Nevertheless, for some applications it is enough to have a weaker result, such as: bounded length of conjugates of elementary generators (see \cite{StepVavDecomp}), bounded length of commutators (see \cite{SivStep},\cite{StepVavLength},\cite{HSVZwidth}), or bounded generation with respect to a larger set of generators. Bounded reduction is a variation of the last property. 

A given Chevalley group $G$ over a given ring is said to have bounded reduction if any element of $G$ can be decomposed as a product of bounded number of elementary generators and one (not necessarily elementary) element from a certain subsystem subgroup. In other words, it means that one can reduce any element to the subsystem subgroup by bounded number of elementary transformations. Without requirement for the number of elementary transformations to be bounded this property is called the surjective stability of the $K_1$-functor. In papers \cite{SteinStability}, \cite{SteinMats}, \cite{PlotProc}, \cite{PlotkinStability}, \cite{PlotkinE7}, and \cite{GvozK1Surj} this problem is considered for rings that satisfy certain conditions on stable rank, absolute stable rank, or other similar conditions. Actually, from the proofs of the theorems in these papers one can recover the bound on the required number of elementary transformations, despite the fact that this bound is not stated in papers explicitly. Therefore, these are results on bounded reduction.

However, conditions on stable rank are still very strong. Even though small Jacobson dimension implies small stable rank, rings with large Jacobson dimension usually fail to have small stable rank. In the present paper, we consider another important class of rings. Namely we take a polynomial ring in arbitrary number of variables with coefficients in a small-dimensional ring. Here we use Krull dimension because the techniques require for dimension to behave well with respect to adding an independent variable. 

Without the bound on the number of elementary transformations similar result for classical groups is known as early surjective stability of the $K_1$-functor. For the special linear group this was proved by Suslin in \cite{SuslinSerreConj}. Similar result for the orthogonal group follows from \cite{SusKopLGPSO}, and for the symplectic group it is proven in \cite{KopejkoSpStab}, see also \cite{GMVLGPSp}, \cite{KopejkoSpDimOne}. Note that if the ring of coefficient is a Dedekind domain or a smooth algebra over a field, then this result for all Chevalley groups follows from the homotopy invariance of the non-stable $K_1$-functor, see \cite{AbeLGP},\cite{StavPolynom},\cite{StavSerreConj}.

In the case of special linear and symplectic groups, there are similar results for Laurent polynomial rings, see \cite{KopejkoLorSL}, \cite{KopejkoLorSp}.

In the paper \cite{VasBounded}, Vaserstein obtained the effective version of the Suslin result, i.e. he proved the bounded reduction for the special linear group over a polynomial ring, and gave this bound explicitly. From this result he deduced that the elementary subgroup of the general linear group over an arbitrary finitely generated commutative ring has Kazhdan's property (T).

In \cite{Shalom}, the basic connection between bounded generation and property (T) has been established and used to estimate the Kazhdan Constants for $\SL_n(\Z)$. Later the bounds for these constants were improved in \cite{KasabovKazdanSL}. In order to deduce property (T) from the Vaserstein's result one needs to refer to \cite{SholomAlgebraisation}. 

In fact, property (T) for Chevalley groups and groups similar to them has already been studied by other methods, see \cite{PropertyTforUniversalLattices},\cite{EJKPropertyT}. However, we believe that the bounded reduction has an independent value, and we aim to study this question for other Chevalley groups. It was noted in the concluding remarks of \cite{VasBounded} that the bounded reduction for the symplectic group follows formally from the case of special linear group. In \cite{GvozBoundedOrt} the bounded reduction for orthogonal groups was established, therefore, closing the problem for classical groups. 

In the present paper we deal with exceptional groups. We prove bounded reduction for the groups of types~$E_6$ and~$E_7$; and we give part of the proof for the group $E_8$.

\smallskip

The main result of the present paper is the following theorem.

{\thm\label{main} Let $C$ be a commutative Noetherian ring and $\dim C=D$. Let $A=C[x_1,\ldots,x_n]$. Let $\Delta\le\Phi$ be one of the following embeddings of root systems
	
	\renewcommand{\theenumi}{\alph{enumi}}
	\begin{enumerate}
	\item $D_5\le E_6$;
	
	\item $E_6\le E_7$;
	\end{enumerate}
	\renewcommand{\theenumi}{\arabic{enumi}}
	
	Assume that
	$$
	D\le \begin{cases}
		3 &\text{ for } D_5\le E_6\tc\\
		4 &\text{ for } E_6\le E_7\tp
	\end{cases}
	$$
	For the case $E_6\le E_7$ assume additionally that $C$ is a Jacobson ring.
	
	\smallskip
	
	 Then every element of the group $G(\Phi,A)$ can be reduced to the subgroup $G(\Delta,A)$ by multiplication from the left by $N$ elementary root elements, where
$$
N=\begin{cases}
	36n^2+(72D+80)n+92 &\text{ for } D_5\le E_6\tc\\
	52n^2+(104D+249)n+244 &\text{ for } E_6\le E_7\tp
\end{cases}
$$}

Therefore, this theorem is an extension of \cite{VasBounded} and \cite{GvozBoundedOrt} to the groups of types~$E_6$ and~$E_7$.  

\smallskip

The paper is organised as follows. In Section \ref{PreliminariesAndNotation}, we give all necessary preliminaries and introduce basic notation. In Section \ref{AbsoluteFlexibleStableRank}, we recall the notion of an absolute flexible stable rank introduced in \cite{GvozBoundedOrt}. In Sections \ref{ReductionToPropositions}, \ref{LowSection} \ref{MonicSection}, and \ref{EliminationSection} we give the proof of the main result.

\section{Preliminaries and notation}
\label{PreliminariesAndNotation}
\subsection*{Rings, ideals and dimensions}

By a ring we always mean associative and commutative ring with unity. 

\smallskip

If $R$ is a ring, then by $R^*$ we denote the set of invertible elements in $R$.

\smallskip

For the elements $r_1$,$\ldots$,$r_k\in R$, we denote by $\<r_1,\ldots,r_k\>$ the ideal in $R$ generated by these elements.
 
\smallskip

In the present paper we use three different notions of a ring dimension.

\begin{itemize}

\item By $\dim R$ we denote Krull dimension of the ring $R$. That is the supremum of the lengths of all chains of prime ideals.

\smallskip

\item By $\Jdim R=\dim\Max(R)$ we denote the dimension of the maximal spectrum $\Max(R)$ of the ring $R$. It is equal to the supremum of the lengths of all chains of such prime ideals that coincide with its Jacobson radical.

\smallskip

\item By $\BSdim R$ we denote the Bass--Serre dimension of a ring $R$. That is the minimal $\delta$ such that $\Max(R)$ is a finite union of irreducible Noetherian subspaces of dimension not greater than $\delta$.

\end{itemize}

Obviously, for a Noetherian ring $R$ we have
$$
\BSdim R\le \Jdim R \le \dim R\tp
$$

\smallskip

The following property of Bass--Serre dimension is well known; see Lemma 4.17 in \cite{BakNonabelian}.

{\lem\label{BSInduction} Let $R$ be a ring with $\BSdim R=d<\infty$. Then it has a finite collection $P_1$,$\ldots$,$P_m$ of maximal ideals, such that for any $s\in R\sm \bigcup_i P_i$ we have $\BSdim R/(s)<d$. In case where $d=0$, this means that $s\in R^*$.}

\subsection*{Chevalley groups}

Let $\Phi$ be a reduced irreducible root system, let $G(\Phi,-)$, be a simply connected Chevalley--Demazure group scheme over $\Z$ of type $\Phi$ (see \cite{ChevalleySemiSimp}), and let $T(\Phi,-)$, be a split maximal torus in it. If $R$ is a commutative ring with unit, the group $G(\Phi,R)$ is called the simply connected Chevalley group of type $\Phi$ over $R$.

For a subset $X$ of a group we denote by $\<X\>$ the subgroup generated by $X$.

To each root $\alpha\in\Phi$ there correspond root unipotent elements $x_\alpha(\xi)$, $\xi\in R$,  elementary with respect to $T$. The group generated by all these elements 
$$
E(\Phi,R)=\left\<x_\alpha(\xi)\colon\alpha\in\Phi\tc\,\xi\in R\right\>\le G(\Phi,R)\tp
$$
is called the elementary subgroup of $G(\Phi,R)$. For any $N\in \N$ we denote by $E(\Phi,R)^{\le N}$ the subset of $E(\Phi,R)$ consisting of elements that can be expressed as the product of no more than $N$ elementary root elements.

Any inclusion of root systems $\Delta\sub\Phi$ induces the homomorphisms $G(\Delta,R)\to G(\Phi,R)$, and $E(\Delta,R)\to E(\Phi,R)$ taking elementary root elements to  elementary root elements.

By $U=U(\Phi,R)$ we denote the subgroup of $E(\Phi,R)$ generated by elementary root elements with positive roots, i.e. the unipotent radical of the standard Borel subgroup.

\subsection*{Basic representations and weight diagrams}

Let us fix an order on $\Phi$, and let $\Phi^+$, $\Phi^-$ and $\Pi=\{\alpha_1\tc\ldots\tc\alpha_l\}$ be the sets of positive, negative, and fundamental roots, respectively. Our numbering of the fundamental roots follows that of~\cite{Bourbaki4-6} and~\cite{atlas}. By $\vpi_1$,$\ldots$,$\vpi_l$ one denotes the corresponding fundamental weights. Let $W=W(\Phi)$ be the Weyl group of the root system $\Phi$.

Recall that an irreducible representation $\pi$ of the complex semisimple Lie algebra $L$ is called {\it basic} (see \cite{Matsumoto69}) if the Weyl group $W=W(\Phi)$ acts transitively on the set $\Lambda^*(\pi)$ of nonzero weights of the representation $\pi$. The set of all weights of the representation $\pi$ we denote by $\Lambda(\pi)$.

In case, where zero is not a weight of the basic representation $\pi$, such representation is called {\it microweight} or {\it minuscule} representations, and the list of these representations is classically known (see \cite{Bourbaki7-9}).

Let $L$ be a simple Lie algebra of type $\Phi$. To each complex representation $\pi$ of the algebra $L$ there corresponds a representation $\pi$ of the Chevalley group $G=G(\Phi,R)$ on the free $R$-module $V_{\pi}=V_{\pi}(R)=V_{\pi}(\Z)\otimes_{\Z}R$ (see \cite{Matsumoto69},\cite{Steinberg85}). If $\pi$ is faithful, then we can identify $G$ with its image under this representation. Thus, for an $g\in G$ and $v\in V$ we write $gv$ for the action of $g$ on $v$.

Decompose the module $V=V_{\pi}$ into the direct sum of its weight submodules
$$
V=V^0\oplus \bigoplus_{\lambda\in\Lambda^*(\pi)} V^{\lambda}\tp
$$
Here all the modules $V^{\lambda}$, $\lambda\in\Lambda^*(\pi)$, are one-dimensional. Matsumoto, \cite[Lemma~2.3]{Matsumoto69}, has shown that there is a special base of weight vectors $e_\lambda\in V^\lambda$, $\lambda\in\Lambda^*(\pi)$, $v^0_{\alpha}\in V^0$, $\alpha\in\Delta(\pi)=\Pi\cap\Lambda^*(\pi)$ such that the action of root unipotents $x_{\alpha}(\xi)$, $\alpha\in\Phi$, $\xi\in R$, is described by the a following simple formulas:
\begin{equation*}
	\begin{array}{lllll}
		\text{i.} & \text{ if } & \lambda\in\Lambda^*(\pi)\tc\, \lambda+\alpha\notin\Lambda(\pi)\tc&\text{ then } & x_{\alpha}(\xi)e^\lambda=e^{\lambda}\tc\\
		\text{ii.} & \text{ if }&\lambda, \lambda+\alpha\in\Lambda^*(\pi)\tc&\text{ then }&x_{\alpha}(\xi)e^\lambda=e^{\lambda}\pm\xi e^{\lambda+\alpha}\tc\\
		\text{iii.}& \text{ if }&\alpha\notin\Lambda^*(\pi)\tc&\text{ then }&x_{\alpha}(\xi)v^0=v^0,\text{ for any }v^0\in V^0\tc\\
		\text{iv.}& \text{ if }&\alpha\in\Lambda^*(\pi)\tc&\text{ then }&x_{\alpha}(\xi)v^{-\alpha}=v^{-\alpha}\pm \xi v^0(\alpha)\pm \xi^2v^{\alpha} \tc\\
		 & & & & x_{\alpha}(\xi)v^0=v^0\pm \xi\alpha_*(v^0)v^\alpha\tp
	\end{array}
\end{equation*}

Weight diagram of the representation $\pi$ is a graph whose vertices correspond to the elements of $\Lambda^*(\pi)\sqcup \Delta(\pi)$; and whose edges labeled be the numbers of fundamental roots shows the action of the corresponding elementary root elements on the weight basis. These diagrams serve as a great visual aid for calculations in Chevalley groups. The details of how to construct and operate with weight diagrams can be found in \cite{atlas} (see also \cite{SteinStability}, \cite{PlotkinE7} and \cite{VavThirdlook}). 

For a fundamental weight $\vpi$, we may consider the basic representation with the highest weight $\vpi$. For simplicity, we call it the representation $\vpi$.

Recall that in the present paper we study bounded reduction for the following embedding of root systems: $D_5\le E_6$, $E_6\le E_7$, and $E_7\le E_8$. We denote by $\Phi$ the bigger system, and by $\Delta$ the smaller one. For the rest of the paper we fix the following representation $\vpi$ of the group $G(\Phi,R)$: 

\renewcommand{\theenumi}{\alph{enumi}}
\begin{enumerate}
	\item $\vpi=\vpi_1$ for $D_5\le E_6$;
	
	\item $\vpi=\vpi_7$ for $E_6\le E_7$;
	
	\item $\vpi=\vpi_8$ for $E_7\le E_8$.
\end{enumerate}
\renewcommand{\theenumi}{\arabic{enumi}}

By $\lambda_1$,$\ldots$,$\lambda_{\dim\vpi}$ we denote the weights of this representation with multiplicities, where numbering of weights for $(E_6,\vpi_1)$ and $(E_7,\vpi_7)$ follows that of \cite{GvozK1Surj} and \cite{PlotkinE7}. We do not need to fix a numbering for $(E_8,\vpi_8)$, but we agree that the highest weight has number $1$, and the lowest weight has number $-1$. 

By $e_1$,$\ldots$,$e_{\dim \vpi}$ we denote the corresponding weight basis. Thus $\lambda_1$ is the highest weight and $e_1$ is the highest weight vector. For $b\in V_{\vpi}$ we denote by $b_i$ or $b_{\lambda_i}$ the corresponding coordinate of $b$ in the basis $e_1$,$\ldots$,$e_{\dim \vpi}$. We will identify $b$ with the column vector with entries $b_i$.

By $V_{\vpi}=V_{\vpi}(R)$ we denote the underlying module of the representation $\vpi$. By $\Um_\vpi R$ we denote the set of unimodular vectors in $V_{\vpi}(R)$, i.e. the set of such vectors $b\in V_{\vpi}(R)$ that the elements $b_i$ generate the unit ideal in $R$. By  $\Eq_{\vpi}$ we denote the set of equations that determine the orbit of the highest weight vector of the representation $\vpi$ (see \cite{Lichtenstein} and \cite{VavThirdlook}). By $\Orb_{\vpi} R$ we denote the set of vectors from $V_{\vpi}(R)$ that satisfy the equations from $\Eq_{\vpi}$. Further set $\Um'_\vpi R=\Um_\vpi R\cap \Orb_{\vpi} R$; and $\Um''_\vpi R=G(\Phi,R)e_1$.

Let $\Sigma_1\le \Phi$ be the set roots that have positive coefficient in simple root $\alpha_i$, $i=1$ for $E_6$ and $i=7$ for $E_7$, and $i=8$ for $E_8$. Therefore, $\Delta\cup \Sigma_1$ is a parabolic set of roots with $\Delta$ being the symmetric part, and $\Sigma_1$ being the special part. Let $U_1$ be the unipotent radical of the corresponding parabolic subgroup, and $U_1^{-}$ be the unipotent radical of the opposite parabolic subgroup.

The following Lemmas can be derived from the proof of Chevalley--Matsumoto decomposition theorem (see \cite{ChevalleySemiSimp},\cite{Matsumoto69},\cite{SteinStability}).

{\lem\label{ChevMatForColumn} Let $b\in \Um' R$ be such that $b_1=1$. Then there exists $u\in U_1^-$ such that $b=ue_1$.}

{\lem\label{ChevMats} Let $g\in G(\Phi,R)$ be such that $(ge_1)_1=1$. Then $$g\in U_1^{-}\cdot G(\Delta,R)\cdot U_1=U_1^{-}U_1G(\Delta,R)\tp$$}

We also need the following lemma.

{\lem\label{Local} If the ring $R$ is semilocal, then $\Um' R=\Um''R$.}
\begin{proof}
	Let $b\in \Um' R$. By Lemma~\ref{ChevMatForColumn}, it is enough to prove that there exists $g\in G(\Phi,R)$ such that $(gb)_1$ is invertible (we may then make it 1 by a toric element). Let $J$ be the Jacobson radical of the ring $R$. The reduction map $E(\Phi,R)\to E(\Phi, R/J)$ is surjective; hence it is enough to find $g\in E(\Phi,R/J)$ such that $(g\ovl{b})_1$ is invertible in $R/J$, where $\ovl{b}$ is the reduction of $b$ (in fact, we have $G=E$ for both $R$ and $R/J$). So we may assume that $J=0$, so $R$ is a product of fields. Moreover, we can look for such $g$ separately for each factor, so we may assume that $R$ is a field. If $b$ has at least one nonzero entry in a position that corresponds to a nonzero weight, then we can take $g$ to be the element of the extended Weyl group that shifts this weight to the highest weight. That concludes the proof for $E_6,\vpi_1$ and $E_7,\vpi_7$ because these representations are minuscule. It remains to consider the case for $E_8,\vpi_8$, where $b\in V^0$. It follows easily from the fact that the lattice $E_8$ is self-dual that we have $x_{\alpha_i}(1)b\notin V^0$ for at least one simple root $\alpha_i$; so the problem is reduced to the previous case.
	\end{proof}

\subsection*{Branching tables}

From the weight diagram it is immediate to read off the branching of the corresponding representation with respect to a subsystem subgroup. In the case where $\Delta=\<\Pi\sm\{\alpha_{h}\}\>$ is the symmetric part of the maximal parabolic subset obtained by dropping the $h$-th fundamental root the procedure is particularly easy. Then the restriction of $\pi$ to $G(\Phi,R)$ looks as follows: one has to cut the diagram of $\pi$ through the bonds with the label $h$.

Given a representation $\pi$ of the group $G(\Phi,R)$ and two fundamental roots $\alpha_{h_1}$,$\alpha_{h_2}\in\Pi$, by "branching table, where vertical lines correspond to cuting through the bonds marked with $h_1$, and horizontal lines correspond to cuting through the bonds marked with $h_2$", we mean the table build as follows: at the upper right corner we write the representation $\pi$; at the remaining cells of the upper row we write the components of restriction of $\pi$ to the group $G(\<\Pi\sm\{\alpha_{h_1}\}\>,-)$; at the remaining cells of the left column we write the components of restriction of $\pi$ to the group $G(\<\Pi\sm\{\alpha_{h_2}\}\>,-)$; and in all the remaining cells we write the intersection of the corresponding restrictions. When this intersection is zero we leave the cell blank; and when the intersection or a component is one-dimensional, we denote it by $\circ$, which refer to the node of the weight diagram. The columns of the table, except the left one, are denoted by bold letters, and the rows except the upper one, are denoted by bold numbers.

Here is the example: the branching table for $(E_6,\vpi_1)$, where vertical lines correspond to cuting through the bonds marked with 1, and horizontal lines correspond to cuting through the bonds marked with 6.

\begin{center}
	\begin{tabular}{lc||c|c|c}
		& & {\bf a}& {\bf b}& {\bf c}\\
		&	$E_6,\vpi_1$ & $\circ$ & $D_5,\vpi_5$ & $D_5,\vpi_1$\\
		\hline\hline
		{\bf 1})&	$\D_5,\vpi_1$ & $\circ$ & $D_4,\vpi_1$ & $\circ$ \\
		\hline
		{\bf 2})&$D_5,\vpi_5$ &	 & $D_4,\vpi_4$ & $D_4,\vpi_3$ \\
		\hline
		{\bf 3})&$\circ$ &  &  & $\circ$
	\end{tabular}
\end{center}

\subsection*{ASR-condition}
Recall that a commutative ring $R$ satisfies the absolute stable rank condition $\ASR_d$ if for any row $(b_1,\ldots,b_d)$ with coordinates in $R$, there exist elements $c_1$,$\ldots$,$c_{d-1}\in R$ such that every maximal ideal of $R$ containing the ideal $\<b_1+c_1b_d,\ldots,b_{d-1}+c_{d-1}b_d\>$ contains already the ideal $\<b_1,\ldots,b_d\>$. This notion was introduced in \cite{EstesOhm} and used in \cite{SteinStability}, \cite{SteinMats} and then in \cite{PlotProc},\cite{PlotkinStability}, \cite{PlotkinE7}, and \cite{GvozK1Surj} to study stability problems.

If we assume that a row $(b_1,\ldots,b_d)$ is unimodular, then the absolute stable rank condition boils down to the usual stable rank condition $\SR_d$ (see \cite{Bass64},\cite{VasersteinStability}).

Absolute stable rank satisfies the usual properties, namely for every ideal $I\unlhd R$ condition $\ASR_d$ for $R$ implies $\ASR_d$ for the quotient $R/I$, and if $d\ge d'$, then $\ASR_{d'}$ implies $\ASR_d$. Finally, it is well known that if the maximal spectrum of $R$ is a Noetherian space of dimension $\Jdim R=d-2$, then both conditions $\ASR_d$ and $\SR_d$ are satisfied (see \cite{EstesOhm},\cite{MagWandarKalVas},\cite{SteinStability}).

\section{Absolute flexible stable rank}
\label{AbsoluteFlexibleStableRank}
In this section, we recall the definition and the basic properties of {\it absolute flexible stable rank} introduced in \cite{GvozBoundedOrt}. Here is the definition.

{\defn A commutative ring $A$ satisfies the absolute flexible stable rank condition $\AFSR_d$ if for any row $(b_1,\ldots,b_d)$ with coordinates in $A$, there exists an element $c_1\in A$ such that for any invertible element $\eps_1\in A^*$, there exists $c_2\in A$ such that for any $\eps_2\in A^*$, $\ldots$, there exists $c_{d-1}\in A$ such that for any $\eps_{d-1}\in A^*$, every maximal ideal of $A$ containing the ideal $\<b_1+\eps_1c_1b_d,\ldots,b_{d-1}+\eps_{d-1}c_{d-1}b_d\>$ contains already the ideal $\<b_1,\ldots,b_d\>$. }

\smallskip

One can think of it as follows. Two players are playing a game. Player 1 chooses a row $(b_1,\ldots,b_d)$ with coordinates in $A$. Then they take turns starting with Player 2. Player 2 in his $i$-th turn chooses element $c_i\in A$; after that Player 1 in his turn chooses invertible element $\eps_i\in A^*$. Player 2 wins if after $d$ turns every maximal ideal of $A$ containing the ideal $\<b_1+\eps_1c_1b_d,\ldots,b_{d-1}+\eps_{d-1}c_{d-1}b_d\>$ contains already the ideal $\<b_1,\ldots,b_d\>$. A commutative ring $A$ satisfies the absolute flexible stable rank condition $\AFSR_d$ if Player 2 has a winning strategy.

The following Lemma shows that the condition $\AFSR_d$ holds for small-dimensional rings. That generalises the result of \cite{EstesOhm}.

{\lem\label{DimmaxImpliesAFSR} (Lemma 3.2 in \cite{GvozBoundedOrt}). Let $A$ be a commutative ring. Assume that $\Max(A)$ is Noetherian and $\Jdim A\le d-2$. Then $A$ satisfies $\AFSR_d$.}

\smallskip

Now the following lemma shows how one can use the $\AFSR$ condition.

{\lem\label{UseAFSR}(Lemma 3.3 in \cite{GvozBoundedOrt}) Let $A$ be a commutative ring, and $S$ be a multiplicative system in $A$. Assume that the localisation $A[S^{-1}]$ satisfies $\AFSR_d$. Then for any row $(b_1,\ldots,b_d)$ with coordinates in $A[S^{-1}]$ and for any $s\in S$, there exist $c_1$,$\ldots$,$c_{d-1}\in sA$ such that every maximal ideal of $A[S^{-1}]$ containing the ideal $\<b_1+c_1b_d,\ldots,b_{d-1}+c_{d-1}b_d\>$ contains already the ideal $\<b_1,\ldots,b_d\>$.}

\section{Reduction of Theorem \ref{main} to propositions \ref{LowDim}, \ref{MakeMonic} and \ref{KillTheVariable}}
\label{ReductionToPropositions}

In this section, we divide the proof of Theorem~\ref{main} into three steps. One of the steps, namely Proposition~\ref{KillTheVariable}, will be formulated and then proved also for the case $E_7\le E_8$, so that if a proof of Propositions \ref{LowDim} and \ref{MakeMonic} are found for this case, it will finish the proof of bounded reduction for Chevalley groups of type $E$.

Recall that by $U_1$ we denote the unipotent radical that correspond to the set $\Sigma_1\le \Phi$, which is the special part of the parabolic subset of roots $\Delta\cup \Sigma_1$. 

Note that
$$
|\Sigma_1|=\begin{cases}
	16 &\text{ for } D_5\le E_6\tc\\
	27 &\text{ for } E_6\le E_7\tp
	\end{cases}
$$

Therefore, Theorem \ref{main} follows trivially from the following result and Lemma~\ref{ChevMats}.

{\thm\label{ColumnReduction} Under the condition of Theorem \ref{main}, for every column $b\in\Um'_{\vpi}A$ there exists a column
$$
b'\in E(\Phi,A)^{\le N}b \tc
$$
where 
$$
N=\begin{cases}
	36n^2+(72D+80)n+60 &\text{ for } D_5\le E_6\tc\\
	52n^2+(104D+249)n+190 &\text{ for } E_6\le E_7\tc
\end{cases}
$$
such that $b'_1=1$.}

\smallskip

Consider the lexicographic order on the monomials in variables $x_1$,$\ldots$,$x_n$. That is the order where $x_1^{k_1}\ldots x_n^{k_n}$ is bigger than $x_1^{l_1}\ldots x_n^{l_n}$ if for some $m$ we have $k_i=l_i$ for $i<m$, and $k_m>l_m$. A polynomial in $A=C[x_1,\ldots,x_n]$ is called lexicographically monic if its leading coefficient in lexicographic order is equal to one.

\smallskip

Further we reduce Theorem \ref{ColumnReduction} to the following three propositions.

{\prop\label{LowDim} Under the condition of Theorem \ref{main}, assuming $n=0$ (i.e. $A=C$), for every column $b\in\Um'_{\vpi}A$ there exists a column
	$$
	b'\in E(\Phi,A)^{\le N}b \tc
	$$
	where 
$$
		N=\begin{cases}
			60 &\text{ for } D_5\le E_6\tc\\
			190 &\text{ for } E_6\le E_7\tp
		\end{cases}
$$
such that $b'_1=1$.}

{\prop\label{MakeMonic} Let $j$ be a number of a vertex on a weight diagram of the representation~$\vpi$. Under the condition of Theorem~\ref{main}, for every column $b\in\Um'_{\vpi}A$ there exists a column
$$
b'\in E(\Phi,A)^{\le N}b\tc
$$
where 
$$
	N=\begin{cases}
		116 &\text{ for } D_5\le E_6\tc\\
		301 &\text{ for } E_6\le E_7\tc
	\end{cases}
$$
such that its entry $b'_{j}$ is lexicographically monic.
}

\smallskip

We state and prove the third proposition also for the case $E_7\le E_8$.

{\prop\label{KillTheVariable}
Let $B$ be a commutative ring such that $\BSdim B=d<\infty$. Let $A=B[y]$, $b\in\Um'_{\vpi}A$ such that its entry $b_{j}$ is monic, where $j=24$ for $E_6$, $j=-1$ for $E_7$ and~$E_8$. Then
$$
E(\Phi,A)^{\le N}b\cap \Um_{\vpi}(B)\ne\emp\tc
$$
where 
$$
	N=\begin{cases}
		72d &\text{ for } D_5\le E_6\tc\\
		104d &\text{ for } E_6\le E_7\tc\\
		291d &\text{ for } E_7\le E_8\tp
	\end{cases}
$$
}

\smallskip

First we need the following lemma.

{\lem\label{Normalisation} Let $f\in C[x_1,\ldots,x_n]$ be a lexicographically monic polynomial. Then there exists an invertible change of variables
$$
x_1,\ldots,x_n \leftrightarrow y_1,\ldots, y_n\tc
$$
such that $f$ becomes monic in $y_n$.}

\begin{proof}
	Take $K>\deg f$. Set $x_i=y_i+y_n^{K^{n-i}}$, $i=1$,$\ldots$,$n-1$, and $x_n=y_n$.
\end{proof}

\smallskip

Now we deduce Theorem \ref{ColumnReduction} from Propositions~\ref{LowDim}, \ref{MakeMonic} and~\ref{KillTheVariable}.

\smallskip

Take $b\in\Um'_{\vpi}A$. By Proposition~\ref{MakeMonic} there exists a column
$$
b'\in E(\Phi,A)^{\le N_1'}b\tc
$$
where 
$$
N_1'=\begin{cases}
	116 &\text{ for } D_5\le E_6\tc\\
	301 &\text{ for } E_6\le E_7\tc
\end{cases}
$$
such that its entry $b'_{j}$ is lexicographically monic, where $j$ is as in Proposition~\ref{KillTheVariable}. Applying Lemma~\ref{Normalisation}, we change variables to $y_1$,$\ldots$,$y_n$ so that $b'_{j}$ is now monic in $y_n$. Now we apply Proposition~\ref{KillTheVariable} to $B=C[y_1,\ldots,y_{n-1}]$. Note that $\BSdim B\le \dim B=D+n-1$. Hence we can obtain a column from
\begin{gather*}
E(\Phi,A)^{\le N_1''}b'\cap \Um_{\vpi}B\le E(\Phi,A)^{\le N_1}b\cap \Um_{\vpi}B\tc
\end{gather*}
where 
$$
N_1''=\begin{cases}
	72(D+n-1) &\text{ for } D_5\le E_6\tc\\
	104(D+n-1) &\text{ for } E_6\le E_7\tc
\end{cases}
$$
 and $N_1=N_1'+N_1''$.

Repeating this argument $n$ times we can obtain a column from
$$
E(\Phi,A)^{\le N_n}b\cap \Um_{\vpi}C\tc
$$
where
\begin{gather*}
N_n=\begin{cases}
	72n(2D+n-1)/2+116n &\text{ for } D_5\le E_6\tc\\
	104n(2D+n-1)/2+301n &\text{ for } E_6\le E_7\tp\\
\end{cases}\tp
\end{gather*}

Now Proposition~\ref{LowDim} implies that there exists
$$
b''\in E(\Phi,A)^{\le N}b\tc
$$
where 
\begin{gather*}
N=\begin{cases}
	72n(2D+n-1)/2+116n+60 &\text{ for } D_5\le E_6\tc\\
	104n(2D+n-1)/2+301n+190 &\text{ for } E_6\le E_7\tc
\end{cases}
\end{gather*}
such that $b''_1=1$.

\section{Bounded reduction for low dimensional rings}
\label{LowSection}

In this Section, we prove Proposition~\ref{LowDim}. Note that the case $E_6$ easily follows from the proof of Lemma 2 in \cite{GvozK1Surj}. Similarly, the proof for the case for the case $E_7$ can be obtain from the proof of the main theorem in \cite{PlotkinE7}. In order to do so, we must estimate how many elementary root element it takes to apply Lemma 2 in \cite{PlotkinE7}. That proof starts with picking an element $e\in E(E_6,R)$ such that $(ae)_{\lambda_1}\equiv 1\mod \u$ and $(ae)_{\lambda_i}\equiv 0\mod \u$ for $i\ne 1$. Note that it is enough to require $(ae)_{\lambda_1}$ to be invertible modulo $\u$ and not necessarily congruent to 1. Therefore, the element $e$ can be taken from $XU_1$ where 
\begin{multline*}
X=\{x_{-\delta_{E_6}}(\xi_1)x_{-\delta_{A_5}}(\xi_2)x_{-\delta_{D_5(6)}}(\xi_3)x_{-\alpha_1}(\xi_4)x_{-\delta_{D_5(1)}}(\xi_5)x_{-\alpha_2-\alpha_3-\alpha_4}(\xi_6)x_{-\alpha_2}(\xi_7)x_{-\alpha_3}(\xi_8)\\x_{-\alpha_4}(\xi_9)x_{-\alpha_5}(\xi_{10})x_{-\alpha_6}(\xi_{11})\colon \xi_i\in R\}\tc
\end{multline*}
where $\delta_{E_6}$ is the maximal root of the system $E_6$; $\delta_{A_5}$ is the maximal root of the system generated by $\alpha_1$,$\alpha_3$,$\alpha_4$,$\alpha_5$ and $\alpha_6$; $\delta_{D_5(6)}$ is the maximal root of the system generated by $\alpha_1$,$\alpha_2$,$\alpha_3$,$\alpha_4$ and $\alpha_5$; $\delta_{D_5(1)}$ is the maximal root of the system generated by $\alpha_2$,$\alpha_3$,$\alpha_4$,$\alpha_5$ and~$\alpha_6$.

The next step in the proof uses the element $e_1\in U_1$, hence $ee_1\in  XU_1\le E^{\le 27}(E_6,R)$. Now it is easy to count that the proof of Lemma 2 takes 67 elementary root elements, and the whole proof of the main theorem in \cite{PlotkinE7} takes 190 elementary root elements.

\section{Obtaining a monic polynomial}
\label{MonicSection}

In this section, we give the proof of Proposition~\ref{MakeMonic}. First we need some preparation.

\smallskip

The following Lemma was proved in \cite{GvozBoundedOrt}. 

{\lem\label{DimOfLocalisation} {\rm (}Lemma 5.2 in \cite{GvozBoundedOrt}{\rm)}. Let $C$ be a Noetherian ring, $A=C[x_1,\ldots,x_n]$. Let $S$ be a multiplicative system of lexicographically monic polynomials in $A$. Then we have $\dim A[S^{-1}]\le \dim C$.}

\smallskip

Now we recall a definition from \cite{VasBounded}.

{\defn\label{mu} Let $A$ be an associative ring with 1, $s$ be a central element of $A$, $l\ge 2$, $v\in A^{l-1}$ (a column over $A$), $u\in\phan^{l-1}A$ (a row over $A$). We define an $l$ by $l$ matrix over $A$ by
$$
\mu(u,s,v)=\begin{pmatrix}
	1_{l-1}+vsu & vs^2\\ -uvu & 1-uvs 
\end{pmatrix}\tp
$$}

This matrix is invertible with $\mu(u,s,v)^{-1}=\mu(u,s,-v)$. If $s\in A^*$, then
$$
\mu(u,s,v)=\begin{pmatrix} 1_{l-1} & 0\\  -u/s & 1 \end{pmatrix}
\begin{pmatrix} 1_{l-1} & vs^2\\  0 & 1 \end{pmatrix}
\begin{pmatrix} 1_{l-1} & 0\\  u/s & 1 \end{pmatrix}
$$ 

The following lemma was proved in \cite{VasBounded}.

{\lem\label{VasersteinMu} {\rm (}Lemma~2.2 in \cite{VasBounded}{\rm )} When $l\ge 3$, the matrix $\mu(u,s,v)$ is a product of $7l-3$ elementary transvections in $\GL(l,R)$.}

\smallskip 

Now let $S$ be a multiplicative system of lexicographically monic polynomials in $A$. It follows from Lemmas~\ref{DimOfLocalisation} and~\ref{DimmaxImpliesAFSR} that the ring $A[S^{-1}]$ satisfies $\AFSR_5$ for the case $E_6$ resp. $\AFSR_6$ for the case $E_7$, and so does any quotient of $A[S^{-1}]$.

{\lem\label{Dr} Let $l\ge 3$, let $\I$ be an ideal in $A$, and suppose that $\dim\Max A/\I[S^{-1}]\le l-2$. Let $b\in V_{D_l,\vpi_1} A$ be such that it becomes unimodular in $ A/\I[S^{-1}]$. Then there exists a column
$$
b'\in E(D_l, A)^{\le 11l-7}b\tc
$$
such that $b_1'$ is congruent to a lexicographically monic polynomial modulo $\I$.}

\begin{proof}
We perform the following steps.

\begin{figure}[h]
	\begin{center}
		\includegraphics[scale=0.4]{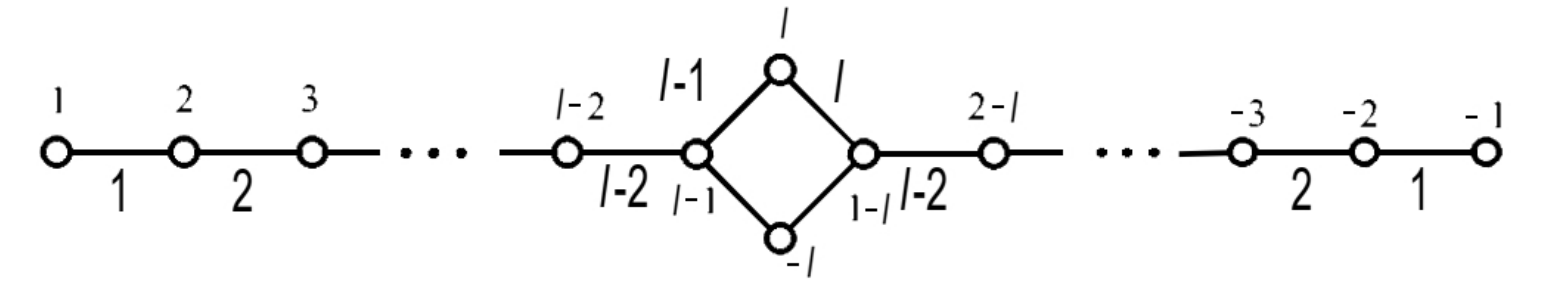}
	\end{center}
	\caption{$(D_l,\vpi_1)$}
	\label{WeightDiagramDl}
\end{figure}

{\bf Step 1.} Make the row $(b_2,\ldots,b_{-1})$ unimodular in $A/\I[S^{-1}]$ by $l-1$ elementary elements.

\smallskip

Let $\A=\I+\<b_{-l},\ldots,b_{-1}\>\unlhd A[S^{-1}]$. Since $A[S^{-1}]/\A$ satisfies $\AFSR_l$ and the row $(b_1,\ldots,b_l)$ is unimodular in $A[S^{-1}]/\A$, it follows from Lemma~\ref{UseAFSR} that there exist $c_2$,$\ldots$,$c_l\in A$ such that the row $(b_2+c_2b_1,\ldots,b_l+c_lb_1)$ is unimodular in $A[S^{-1}]/\A$. Thus by applying the elements $x_{-\alpha_1-\ldots-\alpha_{i-1}}(\pm c_i)$ for $i=2$,$\ldots$,$l$, we make the row $(b_2,\ldots,b_l)$ unimodular in $A[S^{-1}]/\A$ without changing the ideal $\A$. Thus the row $(b_2,\ldots,b_{-1})$ becomes unimodular in $A/\I[S^{-1}]$.

\smallskip

{\bf Step 2.} Make the row $(b_1, b_{-l}\ldots,b_{-1})$ unimodular in $A/\I[S^{-1}]$ by $l-1$ elementary elements.

\smallskip

Since the row  $(b_2,\ldots,b_{-1})$ unimodular in $A/\I[S^{-1}]$, it follows that the ideal generated by $\I$ and $(b_2,\ldots,b_{-1})$ in $A$ contains a lexicographically monic polynomial. So for some $f_2$,$\ldots$,$f_{-1}\in A$, and $f\in \I$ the polynomial
$$
f+\sum_{i=2}^{-1} f_ib_i
$$
is lexicographically monic. Multiplying polynomials $f$ and $f_i$ by a large enough power of $x_1$, we may assume that the polynomial  
$$
b_1+f+\sum_{i=2}^{-1} f_ib_i
$$
is also lexicographically monic.

Let us now apply the elements $x_{\alpha_1+\ldots+\alpha_{i-1}}(f_i)$ for $i=2$,$\ldots$,$l$. Then the ideal generated by $\I$, the new $b_1$, and old $b_{-l}$,$\ldots$,$b_{-1}$ contains a lexicographically monic polynomial. However, these elements do not change the ideal generated by $b_{-l}$,$\ldots$,$b_{-1}$. Hence we actually achieve that the ideal generated by $\I$, and new $b_1$,$b_{-l}$,$\ldots$,$b_{-1}$ contains a lexicographically monic polynomial. Thus the row $(b_1, b_{-l}\ldots,b_{-1})$ becomes unimodular in $A/\I[S^{-1}]$.

\smallskip

{\bf Step 3.} Make the row $(b_1, b_{-l}\ldots,b_{-2})$ unimodular in $A[S^{-1}]$ by $7l-3$ elementary elements.

\smallskip

Let $\A=\I+\<b_{-l},\ldots,b_{-1}\>\unlhd A[S^{-1}]$. Since $b_1$ is invertible in $A[S^{-1}]/\A$, it follows that there exist $\xi_2$,$\ldots$,$\xi_l\in A[S^{-1}]$ such that $b_i-\xi_ib_1\in \A$ for $i=2,\ldots,l$. Let $s$ be a common denominator of $\xi_i$. Set 
$$
g_1=\prod_{2\le i\le l} x_{-\alpha_1-\ldots-\alpha_{i-1}}(\pm\xi_i)\tc
$$
where signs are such that $(g_1b)_i= b_i-\xi_ib_{1}\in \A$ for $1\le i\le l$.

Since $A/\I[S^{-1}]$ satisfies $\AFSR_l$, it follows from Lemma~\ref{UseAFSR} that there exist $c_{-l}$,$\ldots$,$c_{-2}\in s^2A$ such that every maximal ideal of $A/\I[S^{-1}]$ containing the ideal $\<(g_1b)_{-l}+c_{-l}(g_1b)_{-1},\ldots (g_1b)_{-2}+c_{-2}(g_1b)_{-1}\>$ contains already the ideal $\<(g_1b)_{-l},\ldots,(g_1b)_{-1}\>=\A$. Set
$$
g_2=\prod_{2\le i\le l} x_{\alpha_1+\ldots+\alpha_{i-1}}(\pm c_{-i})\tc
$$
where signs are such that $(g_2g_1b)_i= (g_1b)_i+c_ib_{-1}$ for $-l\le i\le-2$.

We claim that the elements $(g_2g_1b)_1$,$(g_2g_1b)_{-l}$,$\ldots$,$(g_2g_1b)_{-2}$ generate the unit ideal in $A/\I[S^{-1}]$. Let us prove that.

Assume that some maximal ideal $\M$ of the ring $A[S^{-1}]$ contains $\I$ and all the elements $(g_2g_1b)_1$,$(g_2g_1b)_{-l}$,$\ldots$,$(g_2g_1b)_{-2}$

Since applying $g_1$ does not change the ideal generated by $b_{-l}$,$\ldots$,$b_{-1}$, by choice of $c_i$ we have $\A \le\M$. Hence $(g_1b)_i\in \M$ for $2\le i\le l$. Thus $b_{1}=(g_2g_1b)_{1}+\sum_{2\le i\le 7}\pm c_{-i}(g_1b)_{i}\in\M$. However, by previous step, $b_{1}$ and $\A$ generate a unit ideal. This is a contradiction.

Since applying $g_1^{-1}$ does not change the ideal generated by elements $b_1$, $b_{-l}$,$\ldots$,$b_{-2}$, we obtain that the elements $(g_1^{-1}g_2g_1b)_i$, where $i=1,-l,\ldots,-2$, generate the unit ideal in $A[S^{-1}]$.

It remains to notice that the element $g_1^{-1}g_2g_1$ is the image of the matrix $\mu(u,s,v)$ for certain $u$ and $v$ under the embedding $G(A_{l-1},A)\to G(D_l,A)$ as a subsystem subgroup. Therefore, by Lemma~\ref{VasersteinMu}, $g_1^{-1}g_2g_1\in E(D_l,A)^{\le 7l-3}$. 

\smallskip

{\bf Step 4.} Make the row $(b_{-l}\ldots,b_{-2})$ unimodular in $A/\I[S^{-1}]$ by $l-1$ elementary elements.

\smallskip

Since $A/\I[S^{-1}]$ satisfies $\AFSR_l$ and the row $(b_1, b_{-l}\ldots,b_{-2})$ is unimodular in $A/\I[S^{-1}]$, it follows from Lemma~\ref{UseAFSR} that there exist $c_{-l}$,$\ldots$,$c_{-2}\in A$ such that the row $(b_{-l}+c_{-l}b_{1},\ldots,b_{-2}+c_{-2}b_{1})$ is unimodular in $A/\I[S^{-1}]$. Thus by applying the elements $x_{\alpha_2+\ldots+\alpha_{i-1}-\delta}(\pm c_{-i})$ for $i=2$,$\ldots$,$l$, we make the row $(b_{-l}\ldots,b_{-2})$ unimodular in $A/\I[S^{-1}]$.

\smallskip

{\bf Step 5.} Make $b_{1}$ congruent to a lexicographically monic polynomial modulo $\I$ by $l-1$ elementary elements.

\smallskip

Since the row  $(b_{-l}\ldots,b_{-2})$ is unimodular in $A/\I[S^{-1}]$, it follows that the ideal generated by $\I$ and $(b_{-l},\ldots,b_{-2})$ in $A$ contains a lexicographically monic polynomial. So for some $f_1$,$f_{-l}$,$\ldots$,$f_{-3}\in A$ and $f\in \I$, the polynomial
$$
f+ \sum_{-l\le i\le -2} f_ib_i
$$
is lexicographically monic. Multiplying polynomials $f_i$ and $f$ by a large enough power of $x_1$, we may assume that the polynomial  
$$
b_{-2}+f+ \sum_{-l\le i\le -2} f_ib_i
$$
is also lexicographically monic.

Now applying the elements $x_{\delta-\alpha_2-\ldots-\alpha_{i-1}}(\pm f_{-i})$ for $i=2$,$\ldots$,$l$, we achieve that $b_{1}$ is congruent to a lexicographically monic polynomial modulo $\I$.
	\end{proof}

{\rem One can notice that the proof above repeats the proof of Proposition~4.2 in~\cite{GvozBoundedOrt}, which on its turn basically repeats the proof of stability theorem for $K_1$-functor given in \cite{SteinStability}.}

Now we prove Proposition~\ref{MakeMonic} for the case $(E_6,\vpi_1)$.  

\begin{proof}
	
Consider the branching table for $(E_6,\vpi_1)$, where vertical lines correspond to cuting through the bonds marked with 1, and horizontal lines correspond to cuting through the bonds marked with 6.

\begin{center}
	\begin{tabular}{lc||c|c|c}
		& & {\bf a}& {\bf b}& {\bf c}\\
		&	$E_6,\vpi_1$ & $\circ$ & $D_5,\vpi_5$ & $D_5,\vpi_1$\\
		\hline\hline
		{\bf 1})&	$\D_5,\vpi_1$ & $\circ$ & $D_4,\vpi_1$ & $\circ$ \\
		\hline
		{\bf 2})&$D_5,\vpi_5$ &	 & $D_4,\vpi_4$ & $D_4,\vpi_3$ \\
		\hline
		{\bf 3})&$\circ$ &  &  & $\circ$
	\end{tabular}
\end{center}

Take $b\in \Um'_{(E_6,\vpi_1)}A$. We need to obtain a lexicographically monic polynomial by 116 elementary elements. Since the Weyl group acts transitively on weights, it does not matter in which position to obtain a lexicographically monic polynomial. Let us make it with $b_1$. We perform the following steps.

\begin{figure}[h]
	\begin{center}
		\includegraphics[scale=0.5]{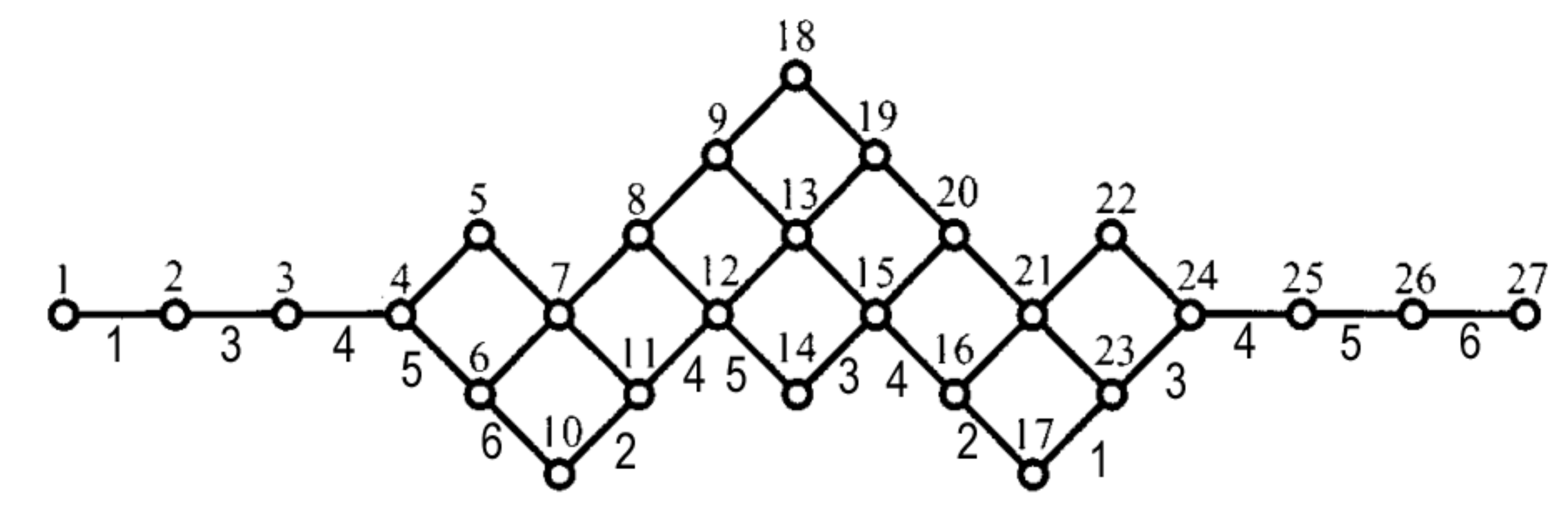}
	\end{center}
	\caption{$(E_6,\vpi_1)$}
	\label{WeightDiagramE6}
\end{figure}

{\bf Step 1.} Make the row that consists of elements in all the cells except {\bf a1} unimodular in $A[S^{-1}]$ by 4 elementary elements.

\smallskip

Let $\A\unlhd A[S^{-1}]$ be the ideal generated by all the elements $b_i$ except for $b_1$,$\ldots$,$b_4$,$b_6$. Since $A[S^{-1}]/\A$ satisfies $\AFSR_5$, and the row $(b_1,\ldots,b_4,b_6)$ is unimodular in $A[S^{-1}]/\A$, it follows from Lemma~\ref{UseAFSR} that there exist $c_2$,$c_3$,$c_4$,$c_6\in A$ such that he row $(b_2+c_2b_1,\ldots,b_6+c_6b_1)$ is unimodular in $A[S^{-1}]/\A$. Thus by applying the elements $x_{\alpha_1}(\pm c_2)$,$\ldots$,$x_{\alpha_1+\alpha_3+\alpha_4+\alpha_5}(\pm c_6)$, we make the row $(b_2,b_3,b_4,b_6)$ unimodular in $A[S^{-1}]/\A$ without changing the ideal $\A$. Thus the row that consists of elements in all the cells except {\bf a1} becomes unimodular in $A[S^{-1}]$.

\smallskip

{\bf Step 2.} Make the row that consists of elements in columns {\bf a} and {\bf c} unimodular in $A[S^{-1}]$ by 16 elementary elements.

\smallskip

Since the row  $(b_2,\ldots,b_{27})$ unimodular in $A[S^{-1}]$, it follows that the ideal generated by $(b_2,\ldots,b_{27})$ in $A$ contains a lexicographically monic polynomial. So for some $f_2$,$\ldots$,$f_{27}\in A$, the polynomial
$$
\sum_{i=2}^{27} f_ib_i
$$
is lexicographically monic. Multiplying polynomials $f_i$ by a large enough power of $x_1$, we may assume that the polynomial  
$$
b_1+\sum_{i=2}^{27} f_ib_i
$$
is also lexicographically monic.

Let us now apply the elements $x_{\lambda_1-\lambda_i}(\pm f_i)$ for all $\lambda_i$ from the column {\bf b}. Then the ideal generated by the new $b_1$ and old $b_{\lambda}$, where $\lambda$ is from the column {\bf c}, contains a lexicographically monic polynomial. However, these elements do not change the ideal generated by $b_{\lambda}$, where $\lambda$ is from the column {\bf c}. Hence we actually achieve that the ideal generated by new $b_1$, and $b_{\lambda}$, where $\lambda$ is from the column {\bf c} contains a lexicographically monic polynomial. Thus the row of elements in columns {\bf a} and {\bf c} becomes unimodular in $A[S^{-1}]$.

\smallskip

{\bf Step 3.} Make the row that consists of elements in cells {\bf a1} and {\bf c1} unimodular in $A[S^{-1}]$ by 48 elementary elements.

\smallskip

Apply Lemma~\ref{Dr} to the column {\bf c} and the ideal generated by $b_1$.

\smallskip

{\bf Step 4.} Make the element $b_1$ lexicographically monic by 48 elementary elements.

\smallskip

Apply Lemma~\ref{Dr} to the row {\bf 1} and the zero ideal.

\end{proof}

{\rem One can notice that the proof above basically repeats the proof of stability theorem for $K_1$-functor given in~\cite{GvozK1Surj}.}

\smallskip

Before we prove Proposition~\ref{MakeMonic} for $(E_7,\vpi_7)$, we need one more Lemma.

{\lem\label{E6} Under the condition of Theorem~\ref{main} in case $E_6\le E_7$. Let $\I$ be an ideal in~$A$. Let $b\in \Orb_{E_6,\vpi_1} A/\I$ be such that it becomes unimodular in $ A/\I[S^{-1}]$. Then there exists a column vector
	$$
	b'\in E(E_6, A)^{\le 91}b\tc
	$$
	such that the row $(b_1',b_{18}')$ is unimodular in $A/\I[S^{-1}]$.}

\begin{proof}
	Let us choose in each irreducible component of $\Max(A/\I[S^{-1}])$ a maximal ideal~$\u_i$, $i\in I$. Next denote by $\tilde{\u_i}$ the preimage of $\u_i$ in $A$. 
	
	For each $i$ choose a maximal ideal $\v_i\in\Max A$ such that it contains $\tilde{\u_i}$, and $b$ is unimodular in $A/\v_i$. Let us show that we can do it. The column $b$ is unimodular in $A/\I[S^{-1}]$; hence the ideal in $A$ generated by $\I$ and entries of $b$ contains a lexicographically monic polynomial $f$. Clearly, the ideal $\tilde{u_i}$ is prime and $f\notin\tilde{\u_i}$. Since $C$ is a Jacobson ring and $A$ is finitely generated over $C$, it follows that $A$ is a Jacobson ring; hence there exists $\v_i\in\Max A$ such that $\tilde{\u_i}\le\v_i$ and $f\notin \v_i$. Then $b$ is unimodular in $A/\v_i$.
	
	Set $I_1=\{i\in I\colon \v_i=\tilde{\u_i}\}$, and $I_2=\{i\in I\colon \v_i\ne\tilde{\u_i}\}$, so $I=I_1\sqcup I_2$.
	
	Now we perform the following steps.
	
	\smallskip
	
	{\bf Step 1.} Achive that $b_1\notin\bigcup_{i\in I}\v_i$ by 11 elementary elements.
	
	 \smallskip
	 
	  In other words, we should make $b_1$ invertible in $A/(\bigcap_{i\in I}\v_i)$. The ring $A/(\bigcap_{i\in I}\v_i)$ is semilocal, so we can do it in the same way we did in Section~\ref{LowSection}.
	  
	 \smallskip
	  
	 {\bf Step 2.} Without changing $b_1$, make the row $(b_2,\ldots,b_{27})$ is unimodular in $ A/\I[S^{-1}]$, and achive that elements $b_5$,$b_7$,$b_8$,$b_9$, $b_{11}$,$\ldots$,$b_{27}$ belong to $\bigcap_{i\in I_1} \v_i$ by element from $U_1$, i.e by 16 elementary elements.
	 
	  \smallskip
	  
	  It follows from Lemma~\ref{ChevMatForColumn} that for some $b'=u_1b$, where $u_1\in U_1$.  we have $b'_{\lambda}\in \bigcap_{i\in I_1} \v_i$ for all $\lambda\ne\lambda_1$. Now let $\A\unlhd A[S^{-1}]$ be the ideal generated by $\I$ and the elements $b'_5$,$b'_7$,$b'_8$,$b'_9$, $b'_{11}$,$\ldots$,$b'_{27}$ except for $b'_1$,$\ldots$,$b'_4$,$b'_6$,$b'_{10}$. Then $\A\le\bigcap_{i\in I_1} \v_i$. Since $A[S^{-1}]/\A$ satisfies $\AFSR_6$, and the row $(b'_1,\ldots,b'_4,b'_6,b'_{10})$ is unimodular in $A[S^{-1}]/\A$, it follows from Lemma~\ref{UseAFSR} that there exist $c_2$,$c_3$,$c_4$,$c_6$,$c_{10}\in A$ such that he row $(b'_2+c_2b'_1,\ldots,b'_{10}+c_{10}b'_1)$ is unimodular in $A[S^{-1}]/\A$. Thus by applying the elements $x_{\alpha_1}(\pm c_2)$,$\ldots$,$x_{\alpha_1+\alpha_3+\ldots+\alpha_6}(\pm c_{10})$, we make the row $(b'_2,b'_3,b'_4,b'_6,b'_{10})$ unimodular in $A[S^{-1}]/\A$ without changing the ideal $\A$. Thus the row $(b'_2,\ldots,b'_{27})$ becomes unimodular in $A/\I[S^{-1}]$.
	  
	  The composition of $u_1$ with the elements as above is then the required element of $U_1$. 
	  
	   \smallskip
	   
	  {\bf Step 3.} Preserving the fact that the image of $b_1$ in $A/\I[S^{-1}]$ does not belong to $\bigcup_{i\in I}\u_i$, make the row that consists of elements in columns {\bf a} and {\bf c} (we use the same branching table as in the proof above) unimodular in $A/\I[S^{-1}]$ by 16 elementary elements.
	   
	    \smallskip
	    
	    Since the row $(b_2,\ldots,b_{27})$ unimodular in $A/\I[S^{-1}]$, it follows that the ideal generated by $\I$ and elements $b_2$,$\ldots$,$b_{27}$ in $A$ contains a lexicographically monic polynomial. So for some $f_2$,$\ldots$,$f_{27}\in A$ and $f\in I$, the polynomial
	    $$
	    f+\sum_{k=2}^{27} f_kb_k
	    $$
	    is lexicographically monic. Clearly for any $i\in I_2$ the ideal $\v_i$ contains some $h_i\in S$. Multiplying polynomials $f$ and $f_k$ by $\prod_{i\in I_2}h_i$ and then by a large enough power of $x_1$, we may assume that the, firstly all the $f_k$ belong to $\bigcap_{i\in I_2}\v_i$, and secondly, that the polynomial  
	    $$
	    b_1+f+\sum_{i=2}^{27} f_ib_i
	    $$
	    is lexicographically monic.
	    
	    Let us now apply the elements $x_{\lambda_1-\lambda_i}(\pm f_i)$ for all $\lambda_i$ from the column {\bf b}. Clearly we preserve the fact that $b_1\notin \bigcup_{i\in I_2} \v_i$; hence we preserve the fact that the image of $b_1$ in $A/\I[S^{-1}]$ does not belong to $\bigcup_{i\in I_2}\u_i$. Further the ideal generated by $\I$, the new $b_1$, and old $b_{\lambda}$, where $\lambda$ is from the column {\bf c}, contains a lexicographically monic polynomial. However, these elements do not change the ideal generated by $b_{\lambda}$, where $\lambda$ is from the column {\bf c}. Hence we actually achieve that the ideal generated by new $b_1$, and $b_{\lambda}$, where $\lambda$ is from the column {\bf c} contains a lexicographically monic polynomial. Thus the row of elements in columns {\bf a} and {\bf c} becomes unimodular in $A[S^{-1}]$. In addition, the ideal generated by $b_{\lambda}$, where $\lambda$ is from the column {\bf c} is contained in $\bigcap_{i\in I_1}\v_i=\bigcap_{i\in I_1}\tilde{\u_i}$, because we make it so at step 2, and it was not changed. Hence the image of new $b_1$ in $A/\I[S^{-1}]$ does not belong to $\bigcup_{i\in I_1}\u_i$.
	    
	    \smallskip
	    
	    {\bf Step 4.} Make the row $(b_1,b_{18})$ unimodular in $A/\I[S^{-1}]$ by 48 elementary elements.
	    
	    \smallskip
	    
	    Since the image of $b_1$ in $A/\I[S^{-1}]$ does not belong to $\bigcup_{i\in I}\u_i$, it follows that $\dim\Max A/(\I+\<b_1\>)[S^{-1}]\le \dim\Max A/\I[S^{-1}]-1\le 3$. Therefore, we can apply Lemma~\ref{Dr} to the column {\bf c} and the ideal $\I+\<b_1\>$.
	\end{proof}

{\rem One can notice that the proof above basically repeats the proof of Lemma~2 in~\cite{PlotkinE7}.}

{\rem The Lemma above is the only place, where we use the assumption that $C$ is a Jacobson ring. It is easy to see that this assumption can be lifted, if we assume that $\dim C\le 3$.}

\smallskip

Now we prove Proposition~\ref{MakeMonic} for the case $(E_7,\vpi_7)$.  

\begin{proof}
	Consider the branching table for $(E_7,\vpi_7)$, where vertical lines correspond to cuting through the bonds marked with 1, and horizontal lines correspond to cuting through the bonds marked with 7.
	
	\begin{center}
		\begin{tabular}{lc||c|c|c}
			& & a& b& c\\
			&	$E_7,\vpi_7$ & $D_6,\vpi_1$ & $D_6,\vpi_6$ & $D_6,\vpi_1$\\
			\hline\hline
			1)&	$\circ$ & $\circ$ & & \\
			\hline
			2)&	$E_6,\vpi_6$ & $D_5,\vpi_1$ & $D_5,\vpi_5$ & $\circ$\\
			\hline
			3)&	$E_6,\vpi_1$ & $\circ$ & $D_5,\vpi_4$ & $D_5,\vpi_1$\\
			\hline
			4)&	$\circ$ & & & $\circ$
		\end{tabular}
	\end{center}
	
	Take $b\in \Um'_{(E_7,\vpi_7)}A$. We need to obtain a lexicographically monic polynomial by 301 elementary elements. Since the Weyl group acts transitively on weights, it does not matter in which position to obtain a lexicographically monic polynomial. Let us make it with $b_1$. We perform the following steps.

	\begin{figure}[h]
		\begin{center}
			\includegraphics[scale=0.5]{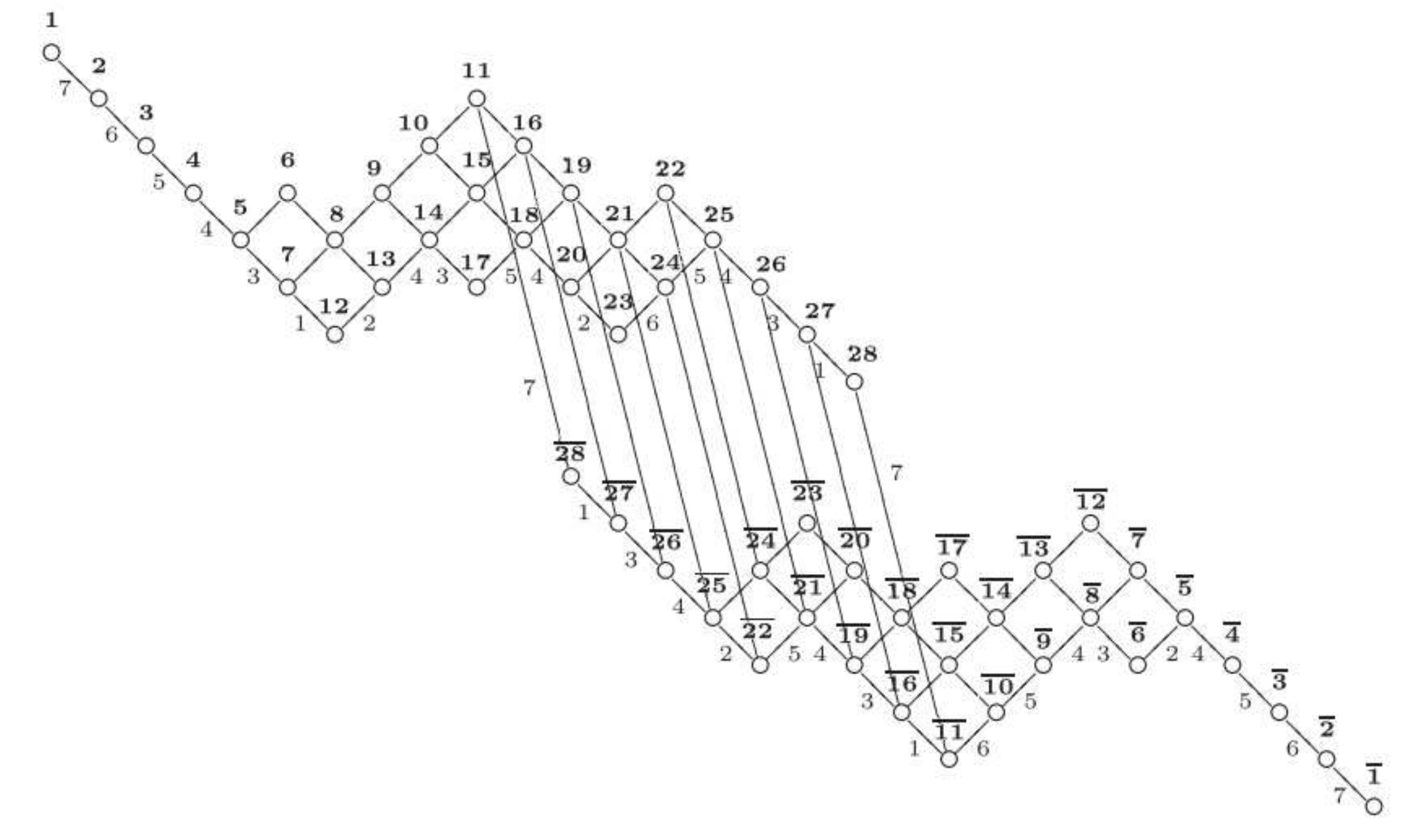}
		\end{center}
		\caption{$(E_7,\vpi_7)$}
		\label{WeightDiagramE76}
	\end{figure}
	
	\smallskip
	
	{\bf Step 1.}  Make the row that consists of elements in all the cells except {\bf a1} unimodular in $A[S^{-1}]$ by 5 elementary elements.
	
	\smallskip
	
	Let $\A\unlhd A[S^{-1}]$ be the ideal generated by all the elements $b_i$ except for $b_1$,$\ldots$,$b_5$,$b_7$. Since $A[S^{-1}]/\A$ satisfies $\AFSR_6$, and the row $(b_1,\ldots,b_5,b_7)$ is unimodular in $A[S^{-1}]/\A$, it follows from Lemma~\ref{UseAFSR} that there exist $c_2$,$\ldots$,$c_5$,$c_7\in A$ such that he row $(b_2+c_2b_1,\ldots,b_{7}+c_{7}b_1)$ is unimodular in $A[S^{-1}]/\A$. Thus by applying the elements $x_{\alpha_7}(\pm c_2)$,$\ldots$,$x_{\alpha_7+\ldots+\alpha_3}(\pm c_{7})$, we make the row $(b_2,\ldots,b_5,b_7)$ unimodular in $A[S^{-1}]/\A$ without changing the ideal $\A$. Thus the row that consists of elements in all the cells except {\bf a1} becomes unimodular in $A[S^{-1}]$.
	
	\smallskip
	
	{\bf Step 2.} Make the row that consists of elements in the rows {\bf 1}, {\bf 3}, and {\bf 4} unimodular in $A[S^{-1}]$ by 27 elementary elements.
	
	\smallskip
	
	Since the row  $(b_2,\ldots,b_{-1})$ unimodular in $A[S^{-1}]$, it follows that the ideal generated by $(b_2,\ldots,b_{-1})$ in $A$ contains a lexicographically monic polynomial. So for some $f_2$,$\ldots$,$f_{-1}\in A$, the polynomial
	$$
	\sum_{i=2}^{-1} f_ib_i
	$$
	is lexicographically monic. Multiplying polynomials $f_i$ by a large enough power of $x_1$, we may assume that the polynomial  
	$$
	b_1+\sum_{i=2}^{-1} f_ib_i
	$$
	is also lexicographically monic.
	
	Let us now apply the elements $x_{\lambda_1-\lambda_i}(\pm f_i)$ for all $\lambda_i$ from the row {\bf 2}. Then the ideal generated by the new $b_1$ and old $b_{\lambda}$, where $\lambda$ is from the rows {\bf 3} and {\bf 4}, contains a lexicographically monic polynomial. However, these elements do not change the ideal generated by $b_{\lambda}$, where $\lambda$ is from the rows {\bf 3} and {\bf 4}. Hence we actually achieve that the ideal generated by new $b_1$, and $b_{\lambda}$, where $\lambda$ is from the rows {\bf 3} and {\bf 4} contains a lexicographically monic polynomial. Thus the row of elements in the rows {\bf 1}, {\bf 3}, and {\bf 4} becomes unimodular in $A[S^{-1}]$.
	
	\smallskip

	{\bf Step 3.} Make the row that consists of elements in the cells {\bf a1}, {\bf b3}, {\bf c3}, and {\bf c4} unimodular in $A[S^{-1}]$ by 5 elementary elements.

	\smallskip
	
	Let $\A\unlhd A[S^{-1}]$ be the ideal generated by all the elements $b_\lambda$ for $\lambda$ in rows {\bf 1}, {\bf 3}, and {\bf 4}, except for $b_{-28}$,$\ldots$,$b_{-23}$. Since $A[S^{-1}]/\A$ satisfies $\AFSR_6$, and the row $(b_{-28},\ldots,b_{-23})$ is unimodular in $A[S^{-1}]/\A$, it follows from Lemma~\ref{UseAFSR} that there exist $c_{-27}$,$\ldots$,$c_{-23}\in A$ such that he row $(b_{-27}+c_{-27}b_{-28},\ldots,b_{-23}+c_{-23}b_{-28})$ is unimodular in $A[S^{-1}]/\A$. Thus by applying the elements $x_{-\alpha_1}(\pm c_{-27})$,$\ldots$,$x_{-\alpha_1-\alpha_3-\ldots-\alpha_6}(\pm c_{-23})$, we make the row $(b_{-27},\ldots,b_{-23})$ unimodular in $A[S^{-1}]/\A$ without changing the ideal $\A$. Thus the row that consists of elements in the cells {\bf a1}, {\bf b3}, {\bf c3}, and {\bf c4} becomes unimodular in $A[S^{-1}]$.
	
	\smallskip
	
	{\bf Step 4.} Make the row that consists of elements in the cells {\bf a1}, {\bf a3}, {\bf c3}, and {\bf c4} unimodular in $A[S^{-1}]$ by 16 elementary elements.
	
	\smallskip
	
	Let $\Gamma$ be the set of weights in cells {\bf a1}, {\bf b3}, {\bf c3}, and {\bf c4} Since the row that consists of elements $\{b_{\lambda}\colon \lambda\in\Gamma\}$ is unimodular in $A[S^{-1}]$, it follows that the ideal generated by $\{b_{\lambda}\colon \lambda\in\Gamma\}$ in $A$ contains a lexicographically monic polynomial. So for some $f_\lambda\in A$, where $\lambda\in \Gamma$ the polynomial
	$$
	\sum_{\lambda\in\Gamma} f_\lambda b_\lambda
	$$
	is lexicographically monic. Multiplying polynomials $f_\lambda$ by a large enough power of $x_1$, we may assume that the polynomial  
	$$
	b_{-28}+\sum_{\lambda\in\Gamma} f_\lambda b_\lambda
	$$
	is also lexicographically monic.
	
	Let us now apply the elements $x_{\lambda_{-28}-\lambda}(\pm f_\lambda)$ for all $\lambda_i$ from the cell {\bf b3}. Then the ideal generated by the new $b_{-28}$ and old $b_{\lambda}$, where $\lambda$ is from the cells {\bf a1}, {\bf c3}, and {\bf c4}, contains a lexicographically monic polynomial. However, these elements do not change the ideal generated by $b_{\lambda}$, where $\lambda$ is from the cells {\bf a1}, {\bf c3}, and {\bf c4}. Hence we actually achieve that the ideal generated by new $b_1$, and $b_{\lambda}$, where $\lambda$ is from the cells {\bf a1}, {\bf c3}, and {\bf c4}, contains a lexicographically monic polynomial. Thus the row of elements in in the cells {\bf a1}, {\bf a3}, {\bf c3}, and {\bf c4} becomes unimodular in $A[S^{-1}]$.
	
	\smallskip
	
	{\bf Step 5.} Make the row that consists of elements in the cells {\bf a1}, {\bf a2}, {\bf a3}, and {\bf c2} unimodular in $A[S^{-1}]$ by 59 elementary elements.
	
	\smallskip
	
	Apply Lemma~\ref{Dr} to the column~{\bf c} and the ideal generated by elements from the column~{\bf a}.
	
	\smallskip

	{\bf Step 6.} Make the row that consists of elements in the cells {\bf a1}, {\bf a2}, and {\bf c2} unimodular in $A[S^{-1}]$ by 39 elementary elements.

	\smallskip
	
	Let $\A\unlhd A[S^{-1}]$ be the ideal generated by elements from column {\bf a}. Since $b_{28}$ is invertible in $A[S^{-1}]/\A$, it follows that there exist $\xi_{-7}$,$\ldots$,$\xi_{-11}\in A[S^{-1}]$ such that $b_i-\xi_ib_{28}\in \A$ for $i=-7$,$\ldots$,$-11$. Let $s$ be a common denominator of $\xi_i$. Set 
	$$
	 g_1=\prod_{-11\le i\le-7} x_{\lambda_i-\lambda_{28}}(\pm\xi_i)\tc
	$$
	where signs are such that $(g_1b)_i= b_i-\xi_ib_{28}\in \A$ for $-11\le i\le-7$.
	
	Since $A[S^{-1}]$ satisfies $\AFSR_6$, it follows from Lemma~\ref{UseAFSR} that there exist $c_{7}$,$\ldots$,$c_{11}\in s^2A$ such that every maximal ideal of $A[S^{-1}]$ containing the ideal $\<(g_1b)_{7}+c_{7}(g_1b)_{-28},\ldots (g_1b)_{11}+c_{11}(g_1b)_{-28}\>$ contains already the ideal $\<(g_1b)_{7},\ldots,(g_1b)_{11},(g_1b)_{-28}\>=\A$. Set
	$$
	g_2=\prod_{7\le i\le 11} x_{\lambda_i-\lambda_{-28}}(\pm c_i)\tc
	$$
	where signs are such that $(g_2g_1b)_i= (g_1b)_i+c_ib_{-28}$ for $7\le i\le 11$.
	
	We claim that the elements $(g_2g_1b)_\lambda$, where $\lambda$ is in the cells {\bf a1}, {\bf a2}, and {\bf c2}, generate the unit ideal in $A[S^{-1}]$ Let us prove that.
	
	Assume that some maximal ideal $\M$ of the ring $A[S^{-1}]$ contains all the elements $(g_2g_1b)_\lambda$, where $\lambda$ is in the cells {\bf a1}, {\bf a2}, and {\bf c2}
	
	Since applying $g_1$ does not change the ideal generated by elements from column {\bf a}, by choice of $c_i$ we have $\A\le\M$. Hence $(g_1b)_i\in \M$ for $-11\le i\le-7$. Thus $b_{28}=(g_2g_1b)_{28}+\sum_{7\le i\le 11}\pm c_i(g_1b)_{-i}\in\M$. However, by previous step, $b_{28}$ and $\A$ generate a unit ideal. This is a contradiction.
	
	Since applying $g_1^{-1}$ does not change the ideal generated by elements from the cells {\bf a1}, {\bf a2}, and {\bf c2}, we obtain that the elements $(g_1^{-1}g_2g_1b)_\lambda$, where $\lambda$ is in the cells {\bf a1}, {\bf a2}, and {\bf c2}, generate the unit ideal in $A[S^{-1}]$
	
	It remains to notice that the element $g_1^{-1}g_2g_1$ is the image of the matrix $\mu(u,s,v)$ for certain $u$ and $v$ under the embedding $G(A_5,A)\to G(E_7,A)$ as a subsystem subgroup. Therefore, by Lemma~\ref{VasersteinMu}, $g_1^{-1}g_2g_1\in E(E_7,R)^{\le 39}$. 
	
	\smallskip
	
	{\bf Step 7.} Make the row that consists of elements in the cells {\bf a1}, and {\bf a2} unimodular in $A[S^{-1}]$ by 91 elementary elements.
	
	\smallskip
	
	It follows from the proof of Lemma~1 in~\cite{PlotkinE7} that the elements in the row {\bf 2}, taken modulo the ideal $\<b_1\>\unlhd A$, form an element of $\Orb_{E_6,\vpi_6} A/\<b_1\>$. Therefore, we can apply Lemma~\ref{E6} to the row {\bf 2} and the ideal $\<b_1\>$. 
	
	\smallskip	
	
{\bf Step 8.} Make the element $b_1$ lexicographically monic by 59 elementary elements.

\smallskip

Apply Lemma~\ref{Dr} to the column {\bf a} and the zero ideal.
	
\end{proof}

\section{Eliminating a variable}
\label{EliminationSection}

In this section, we give the proof of Proposition~\ref{KillTheVariable}. First we need some preparation.

\smallskip

For $\Phi=E_6$,$ E_7$ or $E_8$, let $\Sigma_2\le \Phi$ be the set roots that have positive coefficient in simple root $\alpha_2$, and $\Delta_2$ be the set roots that have zero coefficient in $\alpha_2$. Therefore, $\Delta_2\cup \Sigma_2$ is a parabolic set of roots with $\Delta_2$ being the symmetric part, and $\Sigma_2$ being the special part. Let $U_2$ be the unipotent radical of the corresponding parabolic subgroup, and $U_2^{-}$ be the unipotent radical of the opposite parabolic subgroup. 

Note that
$$
|\Sigma_2|=\begin{cases}
	21 &\text{ for } E_6\tc\\
	42 &\text{ for } E_7\tc\\
	92 &\text{ for } E_8\tp
\end{cases}
$$

Further let $\Lambda$ be the set of weights of the representation $\vpi$. Denote by $\Lambda_i\le \Lambda$ the subset of such weights $\lambda$ that in the decomposition of $\lambda_1-\lambda$ in simple roots the coefficient in $\alpha_2$ is equal to $i$. Therefore, 
$$
\Lambda=\bigcup_{i=0}^{i_{\max}}\Lambda_i\tc
$$   
where
$$
i_{max}=\begin{cases}
	2 &\text{ for } E_6\tc\\
	3 &\text{ for } E_7\tc\\
	6 &\text{ for } E_8\tp
\end{cases}
$$

For an ideal $I\unlhd R$, set
$$
U_2^-(I)=\<x_{\gamma}(\xi)\colon \gamma\text{ has negative coefficient in }\alpha_2\tc\; \xi\in I\>\tp
$$

{\lem\label{ExistsuUnimodCasePart} Let $R$ be a commutative ring. Let $0\le r\le i_{\max}-1$. Let $b$,$b'\in \Orb_{\vpi} R$ be such that for all $0\le i\le r$ and for all $\lambda\in\Lambda_i$ we have $b_{\lambda}=b'_{\lambda}$. Let $I=\<b_{\lambda}-b'_{\lambda}\colon \lambda\in\Lambda\>\unlhd R$. Suppose that the elements $\{b_{\lambda}\colon \lambda\in\Lambda_0\}$ generate the unit ideal. Let $\gamma_1$,$\ldots$,$\gamma_p\in\Phi$ be all the roots with coefficient in $\alpha_2$ being equal to $-(r+1)$. Then there exists an element $u=x_{\gamma_1}(\xi_1)\ldots x_{\gamma_p}(\xi_p)$, where all the $\xi_j$ are in $I$, such that for all $0\le i\le r+1$ and for all $\lambda\in\Lambda_i$ we have $(ub)_{\lambda}=b'_{\lambda}$.}

{\rem In particular, this means that if $r+1$ is bigger than the maximal coefficient in $\alpha_2$, then the assumptions on $b$ and $b'$ implies that $b=b'$.}

\begin{proof}
	Note that such an element $u$ does not change the elements $b_{\lambda}$ for $\lambda\in \bigcup_{i\le r}\Lambda_i$. Therefore, we must ensure the equalities $(ub)_{\lambda}=b'_{\lambda}$ only for $\lambda\in\Lambda_{r+1}$. It is easy to see that this equalities are linear equations in $\xi_j$. We must seek $\xi_j$ in form $\xi_j=\sum_{\lambda} \zeta_{j,\lambda}(b_{\lambda}-b'_{\lambda})$. Therefore, we have a system of linear equations in $\zeta_{j,\lambda}$. For a system of linear equations over a ring, existence of a solution is a local property, see, for example, Proposition~1 in \cite{LinEq}. Hence it is enough to consider the case where $R$ is a local ring.
	
	Note that the system $\Delta_2$ has type $A_{|\Lambda_0|-1}$ and the summand of the representation $\vpi$ that correspond to $\Lambda_0$ is the vector representation of $G(\Delta_2,-)$. Hence in the local case, since the elements $\{b_{\lambda}\colon \lambda\in\Lambda_0\}$ generate the unit ideal, there exists an element $g\in G(\Delta_2,R)\le G(\Phi,R)$ such that $(gb)_{\lambda_1}=1$ and $(gb)_{\lambda}=0$ for all $\lambda\in\Lambda_0\sm \{\lambda_1\}$. The same equalities hold for $gb'$. It follows by Lemma~\ref{ChevMatForColumn} that $gb=u_1e_1$ and $gb'=u_1'e_1$ for some $u_1$,$u_1'\in U_1^-$. Since $(gb)_{\lambda}=(gb')_{\lambda}=0$ for all $\lambda\in\Lambda_0\sm \{\lambda_1\}$, it follows that $u_1$,$u_1'\in U_2^-$. We denote by $U_2^{-(j)}$ the subgroup of $U_2^-$ generated by root subgroups for all roots that have coefficient in $\alpha_2$ less or equal than $-(j+1)$. Similarly, we define the subgroup $U_2^{-(j)}(I)$. Since $g\in G(\Delta_2,R)$, the assumptions on $b$ and $b'$ imply that $(gb)_{\lambda}=(gb')_{\lambda}$ for all $\lambda\in \bigcup_{i\le r}\Lambda_i$; hence we have $u_1\equiv u_1'\mod U_2^{-(r)}$. Moreover, it is easy to see that the ideal generated by elements $(gb)_{\lambda}-(gb')_{\lambda}$ for all $\lambda\in\Lambda$ is equal to $I$. Therefore, $u_1\equiv u_1'\mod U_2^{-(r)}(I)$. Set $\tilde u=u_1'u_1^{-1}\in  U_2^{-(r)}(I)$. Then we have
	$$
	b'=g^{-1}gb'=g^{-1}u_1'e_1=g^{-1}\tilde{u}u_1e_1=g^{-1}\tilde{u}gb\tp
	$$
	Since $g\in G(\Delta_2,R)$, it follows that $g^{-1}\tilde{u}g\in  U_2^{-(r)}(I)$. Therefore, $g^{-1}\tilde{u}g=u\hat{u}$, where $u=x_{\gamma_1}(\xi_1)\ldots x_{\gamma_p}(\xi_p)$, where $\xi_j\in I$, and $\hat{u}\in U_2^{-(r+1)}$. Then we have $(ub)_{\lambda}=(g^{-1}\tilde{u}gb)_{\lambda} =b'_{\lambda}$ for $\lambda\in\Lambda_{r+1}$.
	\end{proof}

{\lem\label{ExistsuUnimodCase} Let $R$ be a commutative ring. Let $b$,$b'\in \Orb_{\vpi} R$ be such that for all $\lambda\in\Lambda_0$ we have $b_{\lambda}=b'_{\lambda}$. Let $I=\<b_{\lambda}-b'_{\lambda}\colon \lambda\in\Lambda\>\unlhd R$. Suppose that the elements $\{b_{\lambda}\colon \lambda\in\Lambda_0\}$ generate the unit ideal. Then there exists an element $u\in U_2^-(I)$ such that $ub=b'$.}

\begin{proof}
	Follows from Lemma~\ref{ExistsuUnimodCasePart} by induction.
	\end{proof}

{\lem\label{ExistsuIfDivisible} Let $R$ be a commutative Noetherian ring, $s\in R$. Then there exists $k\in \N$ such that for any $b$,$b'\in \Orb_{\vpi} R$ that satisfy the following conditions:
	\begin{itemize}
		 \item for all $\lambda\in\Lambda_0$ we have $b_{\lambda}=b'_{\lambda}$,
		 \item $s\in \<b_\lambda\colon \lambda\in\Lambda_0\>\unlhd R$,
		 \item  $b_{\lambda}-b'_{\lambda}$ is divisible by $s^k$ for all $\lambda\in\Lambda$,
\end{itemize}
there exists an element $u\in U_2^-$ such that $ub=b'$.
}

\begin{proof}
	The annihilators of the elements $s^i$, $l\in\N$ form a ascending chain
$$
\Ann s\le \Ann s^2\le \ldots\tp
$$  
Since the ring $R$ is Noetherian, it follows that for some $l\in \N$, we have $\Ann s^{l+m}=\Ann s^l$ for any $m\in\N$.

Now consider the ring $\Z[\{\tilde{b}_{\lambda}\colon \lambda\in\Lambda\}][\{\tilde{b'}_{\lambda}\colon \lambda\in\Lambda\}][\{\tilde{a}_{\lambda}\colon \lambda\in\Lambda_0\}]$ of polynomials over $\Z$ in $2|\Lambda|+|\Lambda_0|$ variables. Set
$$
\tilde{R}=\Z[\{\tilde{b}\colon \lambda\in\Lambda\}][\{\tilde{b'}\colon \lambda\in\Lambda\}][\{\tilde{a}\colon \lambda\in\Lambda_0\}]/\mathfrak{I}\tc
$$
where the ideal $\mathfrak{I}$ is generated by the following elements: equations form $\Eq_{\vpi}$ for the column $\tilde{b}$, equations form $\Eq_{\vpi}$ for the column $\tilde{b'}$, elements $\tilde{b}_{\lambda}-\tilde{b'}_{\lambda}$ for all $\lambda\in\Lambda_0$. Then set
$$
\tilde{s}=\sum_{\lambda\in\Lambda_0}\tilde{a}_{\lambda}\tilde{b}_\lambda\tp
$$

It follows by Lemma~\ref{ExistsuUnimodCase} that over the ring $\tilde{R}[\tilde{s}^{-1}]$ there exists an element $\tilde{u}\in U_2^-$ such that $\tilde{u}\tilde{b}=\tilde{b'}$ in $\tilde{R}[\tilde{s}^{-1}]$. Let $\tilde{u}=x_{\gamma_1}(\tilde{\xi_1})\ldots x_{\gamma_q}(\tilde{\xi_q})$, where $\gamma_j$ are roots that have negative coefficient in $\alpha_2$. Moreover, we can choose $\tilde{\xi_j}$ to be in the ideal generated by $\tilde{b}_{\lambda}-\tilde{b'}_{\lambda}$ in $\tilde{R}[\tilde{s}^{-1}]$. Let $\tilde{k}\in\N$ be such that for all $i$ elements $\tilde{s}^{\tilde{k}}\tilde{\xi_j}$ belong to the ideal generated by $\tilde{b}_{\lambda}-\tilde{b'}_{\lambda}$ in $\tilde{R}$, i.e
$$
\tilde{\xi_j}=\tilde{s}^{-\tilde{k}}\sum_{\lambda\in\Lambda}\tilde{\zeta}_{j,\lambda}(\tilde{b}_{\lambda}-\tilde{b'}_{\lambda})\tc\qquad\tilde{\zeta}_{j,\lambda}\in\tilde{R}\tp
$$

We claim that we can take $k=\tilde{k}+l$. Let $b$,$b'\in\Orb_{\vpi} R$ satisfy the conditions. Then there exists a ring homomorphism $\map{\ph}{\tilde{R}}{R}$ that maps $\tilde{b}$ to $b$, $\tilde{b'}$ to $b'$, and $\tilde{s}$ to $s$. Set $\zeta_{j,\lambda}=\ph(\tilde{\zeta}_{j,\lambda})$. Let $b_{\lambda}-b'_{\lambda}=s^kc_{\lambda}$ for all $\lambda\in\Lambda$. We claim that we can take $u=x_{\gamma_1}(\xi_1)\ldots x_{\gamma_q}(\xi_q)$, where 
$$
\xi_j=s^l\sum_{\lambda\in\Lambda}\zeta_{j,\lambda}c_{\lambda}\tp
$$
It is easy to see that the homomorphism $\tilde{R}[\tilde{s}^{-1}]\to R[s^{-1}]$ induced by $\ph$ send $\tilde{u}$ to $u$. Therefore, $ub=b'$ over $R[s^{-1}]$, i.e. for any $\lambda\in\Lambda$ we have $(ub)_\lambda-b'_{\lambda}\in\Ann s^m$ for some $m$. On the other hand, since all the $\xi_j$ are divisible by $s^l$, it follows that $(ub)_\lambda-b_{\lambda}$ are divisible by $s^l$; hence $(ub)_\lambda-b'_{\lambda}$ are divisible by $s^l$. Let $(ub)_\lambda-b'_{\lambda}=s^l\theta_\lambda$. Then $s^{m+l}\theta_\lambda=s^m((ub)_\lambda-b'_{\lambda})=0$, i.e $\theta_\lambda\in \Ann s^{m+l}=\Ann s^l$. Therefore, $(ub)_\lambda-b'_{\lambda}=s^l\theta_\lambda=0$, i.e $ub=b'$.
	\end{proof}

Set $\Lambda_0'=\Lambda_0\sm\{\nu\}$, where $\nu$ is the lowest weight in $\Lambda_0$.

{\lem\label{AddToY} Let $B$ be a commutative Noetherian ring, $A=B[y]$. Let $b=b(y)\in \Um' A$, and $s\in B\cup \<b(y)_{\lambda}\colon \lambda\in\Lambda_0'\>$. Then there exists $m\in\N$ such that
$$
b(y+s^mz)\in E(\Phi,A[z])^{\le N}b(y)\tc
$$
where 
\begin{gather*}
	N=\begin{cases}
		65 &\text{ for } D_5\le E_6\tc\\
		94 &\text{ for } E_6\le E_7\tc\\
		152 &\text{ for } E_7\le E_8\tp
	\end{cases}
\end{gather*}}

\begin{proof}
	Take $k$ from Lemma~\ref{ExistsuIfDivisible} (for $R=A[z]$) and set $m=k+2$. Let $b=(\begin{smallmatrix}
		b^0\\ b^1
	\end{smallmatrix})$, where $b_0$ is a column with entries $b_{\lambda}$ for $\lambda\in\Lambda_0$ and $b_1$ is a column with entries $b_{\lambda}$ for $\lambda\notin\Lambda_0$. Recall that the system $\Delta_2$ has type $A_{|\Lambda_0|-1}$ and the summand of the representation $\vpi$ that correspond to $\Lambda_0$ is the vector representation of $G(\Delta_2,-)$. Therefore, it follows from Corollary 2.4 of \cite{VasBounded} that there exists an element $g(z)\in E(\Delta_2,A[z])^{\le 8|\Lambda_0|-4}$ such that $g(z)b^0(y)=b^0(y+s^2z)$. Moreover, it follows from the proof that $g$ is congruent to the identity element modulo $z$, see the proof of Lemma~6.5 in \cite{GvozBoundedOrt} for details.

	Therefore, for some $\tilde{b}^1\in A[z]^{|\Lambda\sm\Lambda_0|}$ we have
	$$
	b'=\begin{pmatrix} b^0(y+s^mz) \\ \tilde{b}^1\end{pmatrix}=g(s^kz)b\in E(\Phi,A[z])^{\le 8|\Lambda_0|-4}b\tp
	$$
	In addition, $\tilde{b}^1$ is congruent to $b^1(y)$, and hence to $b^1(y+s^mz)$, modulo $s^k$.	Now applying Lemma~\ref{ExistsuIfDivisible} to vectors $b'$ and $b(y+s^mz)$, we obtain that
	$$
	b(y+s^mz)\in E(\Phi,A[z])^{\le|\Sigma_2|}b'\sub E(\Phi,A[z])^{\le 8|\Lambda_0|-4+|\Sigma_2|}b\tp
	$$
	Here we used that $s\in B$, so the shift of the variable does not change the fact that  $s\in \<b_{\lambda}\colon \lambda\in\Lambda_0\>$.
	\end{proof}

{\lem \label{Weyl} There is an element $w\in W(E_8)$ such that $w(\alpha_2)=\alpha_8$, $w(\alpha_4)=\alpha_7$, $w(\alpha_5)=\alpha_6$, and $w(-\delta_{A_8})=\delta$, where $\delta_{A_8}$ is the maximal root of the subsystem generated by $\alpha_1$,$\alpha_3$,$\ldots$,$\alpha_8$,$\delta$.}

\begin{proof}
	Let $\delta_{A_7}$ be is the maximal root of the subsystem generated by $\alpha_1$,$\alpha_3$,$\ldots$,$\alpha_8$. Let us show that $w'\in W(E_8)$ such that $w(\alpha_2)=\alpha_8$, $w(\alpha_4)=\alpha_7$, $w(\alpha_5)=\alpha_6$, and $w(-\delta_{A_7})=\delta$, where $\delta_{A_8}$. Note that all the roots in the condition belong to the subsystem of type $D_8$ generated by $\alpha_2$,$\alpha_3$,$\ldots$,$\alpha_8$,$\delta$. We can realise this $D_8$ in the Euclidean space with the orthonormal basis $e_1$,$\ldots$,$e_8$ so that $\delta=e_1-e_2$, $\alpha_8=e_2-e_3$, $\alpha_7=e_3-e_4$, $\alpha_6=e_4-e_5$, $\alpha_5=e_5-e_6$, $\alpha_4=e_6-e_7$, $\alpha_2=e_7-e_8$, $\alpha_3=e_7+e_8$, $\delta_{A_7}=e_1+e_8$. An element form $W(D_8)$ can perform any permutation of $e_i$ and in addition replace any even number of $e_i$ with $-e_i$. So we can take $w'\in W(D_8)$ to be the element such that $w'(e_1)=e_1$, $w'(e_2)=e_8$, $w'(e_3)=e_7$, $w(e_4)=e_6$, $w(e_5)=-e_5$, $w(e_6)=-e_4$, $w(e_7)=-e_3$, $w(e_8)=-e_2$.
	
	Now we can take $w=w'w_{\alpha_1}$, where $w_{\alpha_1}$ os the reflection with respect to $\alpha_1$.
	\end{proof}

{\lem \label{PrepareToAddToY} Let $B$ be a commutative ring, $P_1$,$\ldots$,$P_m$ be distinct maximal ideals in $B$, $A=B[y]$, $b=b(y)\in \Um'_{\vpi} A$ such that $b_{j}$ is monic, where $j=24$ for $E_6$, $j=-1$ for $E_7$ and $E_8$. Then there exists a column vector
	$$
	b^{(1)}\in E(\Phi,A)^{\le N}b
	$$
	such that $b^{(1)}_{j}$ is monic and
	$$
	\left(\<b^{(1)}_\lambda\colon \lambda\in\Lambda_0'\>\cap B\right)\sm \bigcup_{i=1}^m P_i\ne \emp\tc
	$$
	where
	\begin{gather*}
		N=\begin{cases}
			7 &\text{ for } D_5\le E_6\tc\\
			10 &\text{ for } E_6\le E_7\tc\\
			139 &\text{ for } E_7\le E_8\tp
		\end{cases}
	\end{gather*}
}

\begin{proof}
	Set 
	$$
	R=B/\left(\bigcap_{i=1}^m P_i\right)=\prod_{i=1}^m B/P_i\tp
	$$
	
	First we show that the last condition on $b^{(1)}$ holds if $b^{(1)}_1$ is monic and the elements $\{b^{(1)}_\lambda\colon\lambda\in\Lambda_0'\}$ generate the unit ideal in $R[y]$.
		
	Let $c_\lambda\in A$, where $\lambda\in\Lambda_0'$, be such that $\sum_{\lambda\in\Lambda_0'}c_\lambda b^{(1)}_\lambda\equiv 1\mod P_i$ for every $i$.
	
	Set $f=\sum_{\lambda\in\Lambda_0'\sm\{\lambda_1\}}c_\lambda b^{(1)}_\lambda$. Then $b^{(1)}_1$ and $f$ are coprime in $B/P_i[y]$ for every $i$.
	
	Since $b^{(1)}_1$ is monic, it follows that the resultant $\res(b^{(1)}_1,f)$ modulo $P_i$ is equal to the resultant of $b^{(1)}$ taken modulo $P_i$ and $f$ taken modulo $P_i$ (even if $f$ modulo $P_i$ has smaller degree).
	
	Therefore, we have 
	$$
	\res(b^{(1)}_1,f)\in \left(\<b^{(1)}_\lambda\colon \lambda\in\Lambda_0'\>\cap B\right)\sm \bigcup_{i=1}^m P_i\tp
	$$
	
	Thus it remains to prove that a given column $b\in \Um'_{\vpi} A$, with $b_{j}$ being monic, can be transformed by $N$ elementary elements so that $b_1$ becomes monic, $b_{j}$ remains monic, and the new elements $\{b_\lambda\colon \lambda\in\Lambda_0'\}$ generate the unit ideal in $R[y]$.
	
	\smallskip
	
	{\bf Proof for $(E_6,\vpi_1)$.} Here we perform the following steps.
	
	\begin{figure}[h]
		\begin{center}
			\includegraphics[scale=0.5]{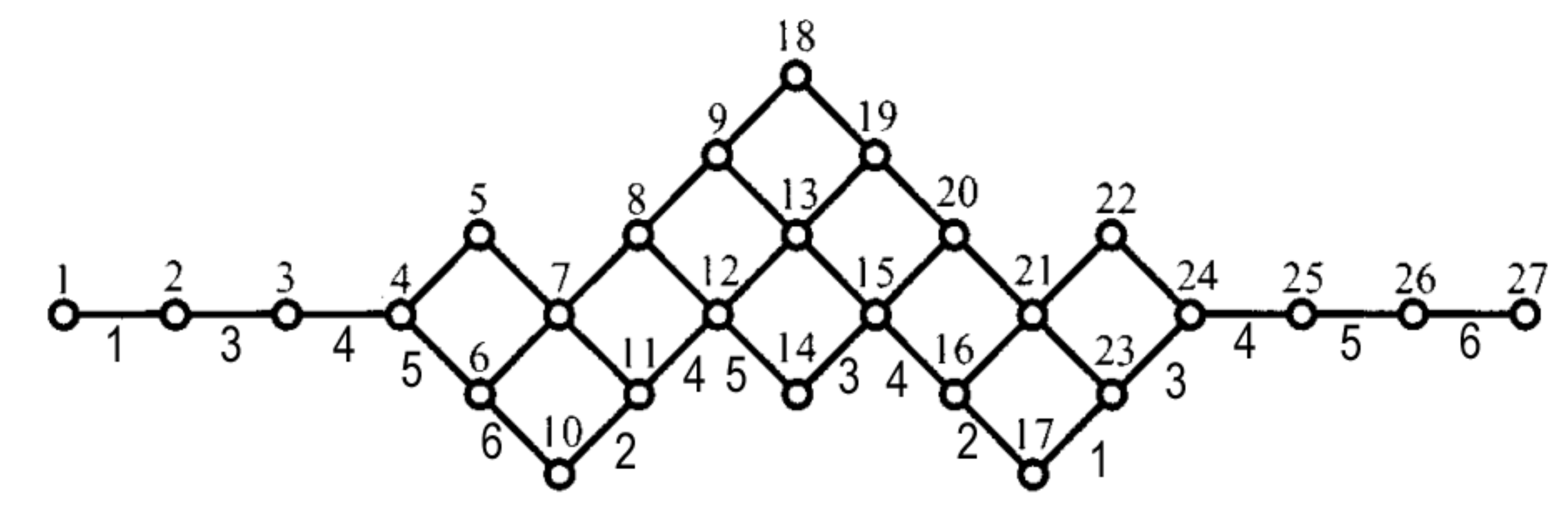}
		\end{center}
		\caption{$(E_6,\vpi_1)$}
		\label{WeightDiagramE6bis}
	\end{figure}
	
	{\bf Step 1.} Make the polynomial $b_3$ monic and the row $(b_1,\ldots,b_{16},b_{18},\ldots b_{22},b_{24})$ unimodular in $R[y]$ by the element $x_{\delta}(\xi)$.
	
	\smallskip
	
	Since $R$ is a product of fields and $b_{24}$ is monic, it follows that the ring $R[y]/\<b_{24}\>$ is semilocal; hence it is easy to see that there exists $\tilde{\xi}$ such that the row $((x_{\delta}(\tilde{\xi})b)_1,\ldots,(x_{\delta}(\tilde{\xi})b)_{16},(x_{\delta}(\tilde{\xi})b)_{18},\ldots (x_{\delta}(\tilde{\xi})b)_{22},(x_{\delta}(\tilde{\xi})b)_{24})$ is unimodular in $R[y]$. Therefore, if we take 
	$$
	\xi=\tilde{\xi}+y^Kb_{24}
	$$
	for some $K\in\N$, then we guarantee that the row $(b_1,\ldots,b_{16},b_{18},\ldots b_{22},b_{24})$ becomes unimodular in $R[y]$. It remains to notice that if $K$ is large enough, then we also make $b_3$ monic.
	
	\smallskip
	
		{\bf Step 2.} Make the polynomial $b_2$ monic and the row $(b_1,\ldots b_{21},b_{23})$ unimodular in $R[y]$ by the element $x_{\alpha_3}(\xi)$.
	
	\smallskip
	
	This is done similarly to step 1.
	
	\smallskip
	
	{\bf Step 3.} Make the polynomial $b_1$ monic and the row $(b_1,\ldots b_{17})$ unimodular in $R[y]$ by the element $x_{\alpha_1}(\xi)$.
	
	\smallskip
	
	This is done similarly to step 1.
	
	\smallskip
	
	{\bf Step 4.} Make the row $(b_1,b_2,b_3,b_4, b_6)$ unimodular in $R[y]$ by the element $x_{\alpha_2}(\xi_4)x_{\delta_{D_4}}(\xi_3)x_{\alpha_6}(\xi_2)x_{\delta_{D_5}}(\xi_1)$, where $\delta_{D_5}$ is the maximal root of the subsystem generated by $\alpha_2,\ldots\alpha_6$, and $\delta_{D_4}$ is the maximal root of the subsystem generated by $\alpha_2,\ldots\alpha_5$.
	
	\smallskip
	
	Existence of such $\xi_1$,$\ldots$,$\xi_4$ follows easily from the fact that $R[y]/\<b_1\>$ is semilocal.
	
	\smallskip
	
	Note that neither of steps change $b_{24}$; hence it remains monic. Also step 4 does not change $b_1$; hence it remains monic after being made so in step 3.
	
	\medskip
	
	{\bf Proof for $(E_7,\vpi_7)$.}
	
	\smallskip
	
	Consider the branching table for $(E_7,\vpi_7)$, where vertical lines correspond to cuting through the bonds marked with 1, and horizontal lines correspond to cuting through the bonds marked with 7.
	
	\begin{center}
	\begin{tabular}{lc||c|c|c}
		 & & a& b& c\\
	&	$E_7,\vpi_7$ & $D_6,\vpi_1$ & $D_6,\vpi_6$ & $D_6,\vpi_1$\\
		\hline\hline
	1)&	$\circ$ & $\circ$ & & \\
		\hline
	2)&	$E_6,\vpi_6$ & $D_5,\vpi_1$ & $D_5,\vpi_5$ & $\circ$\\
		\hline
	3)&	$E_6,\vpi_1$ & $\circ$ & $D_5,\vpi_4$ & $D_5,\vpi_1$\\
		\hline
	4)&	$\circ$ & & & $\circ$
		\end{tabular}
	\end{center}

	Now we perform the following steps, which are similar to those for $E_6$.
	
	\smallskip
	
	{\bf Step 1.} Make the polynomial in the cell {\bf a3} monic and the row that consists of elements in the cells {\bf a1}, {\bf a2}, {\bf a3}, {\bf b2}, {\bf b3}, and {\bf c4}  unimodular in $R[y]$ by the element $x_{\delta}(\xi)$.
	
	\smallskip
	
	{\bf Step 2.} Make the polynomial $b_2$ (highest weight in the cell {\bf a2}) monic and the row that consists of elements in the cells {\bf a1}, {\bf a2}, {\bf a3}, {\bf b2}, {\bf c2}, the upper half of the cell {\bf b3} with respect to cutting through the bonds marked with 6, and the element that correspond to the highest weight in the cell {\bf c3} unimodular in $R[y]$ by the element $x_{\delta_{D_6}}(\xi)$, where $\delta_{D_6}$ is the maximal root of the subsystem generated by $\alpha_2,\ldots\alpha_7$.
	
	\smallskip
	
	{\bf Step 3.} Make the polynomial $b_1$ monic and the row that consists of elements in the cells {\bf a1}, {\bf a2}, {\bf b2}, and {\bf c2} unimodular in $R[y]$ by the element $x_{\alpha_7}(\xi)$.
	
	\smallskip
	
	{\bf Step 4.} Make the row $(b_1,b_2,b_3,b_4, b_5, b_7)$ unimodular in $R[y]$ by the element $x_{\alpha_2}(\xi_7)x_{\alpha_2+\alpha_3+\alpha_4}(\xi_6)x_{\delta_{D_5(1)}}(\xi_5)x_{\alpha_1}(\xi_4)x_{\delta_{D_5(6)}}(\xi_3)x_{\delta_{A_5}}(\xi_2)x_{\delta_{E_6}}(\xi_1)$, where $\delta_{E_6}$ is the maximal root of the system generated by $\alpha_1$,$\alpha_2$,$\alpha_3$,$\alpha_4$,$\alpha_5$ and $\alpha_6$; $\delta_{A_5}$ is the maximal root of the system generated by $\alpha_1$,$\alpha_3$,$\alpha_4$,$\alpha_5$ and $\alpha_6$; $\delta_{D_5(6)}$ is the maximal root of the system generated by $\alpha_1$,$\alpha_2$,$\alpha_3$,$\alpha_4$ and $\alpha_5$; $\delta_{D_5(1)}$ is the maximal root of the system generated by $\alpha_2$,$\alpha_3$,$\alpha_4$,$\alpha_5$ and $\alpha_6$.
	
	\medskip
	
	{\bf Proof for $(E_8,\vpi_8)$.} Here we perform the following steps.
	
	\smallskip
	
	{\bf Step 0.} Make the row $(b_1,b_{-1})$ unimodular in $R[y]$ by the element $u\in U$.
	
	\smallskip
	
Since $R[y]/\<b_{-1}\>$ is semilocal, by Lemma~\ref{Local} there exists $g\in G(E_8,R[y]/\<b_{-1}\>)$ such that $gb=e_1$ in $R[y]/\<b_{-1}\>$. By Theorem 1.1 in~\cite{SmolSuryVavBis}, we have $g=hu_1vu$, where $h\in T$, $u,u_1\in U$, and $v\in U^-$. Therefore, we have $ub=v^{-1}u_1^{-1}h^{-1}e_1$ in $R[y]/\<b_{-1}\>$; hence $(ub)_1$ is invertible in $R[y]/\<b_{-1}\>$. Clearly $u$ can be lifted to the element of $U(\Phi, B[y])$. Note that $(ub)_{-1}=b_{-1}$; hence the row $((ub)_1,(ub)_{-1})$ is unimodular in $R[y]$.

\smallskip
	
	Now consider the subsystem $D_8\le E_8$ generated by $\alpha_2$,$\alpha_3$,$\ldots$,$\alpha_8$,$\delta$. If we restrict our representation to the group $G(D_8,-)$ one of the summands will be the representation $(D_8,\vpi_8)$. Take $w\in W(E_8)$ from Lemma~\ref{Weyl}. If we move our subsystem $D_8$ with element $w$, then the highest weight of the representation $(D_8,\vpi_8)$ become the highest weight of the entire $(E_8,\vpi_8)$. In addition, three weights next to it become weights from $\Lambda_0'$ (Figure~\ref{PartWeightDiagramE8}). It is clear that lowest weight of $(D_8,\vpi_8)$ become the lowest weight of $(E_8,\vpi_8)$.
	
	\begin{figure}[h]
		\begin{center}
			\includegraphics[scale=1]{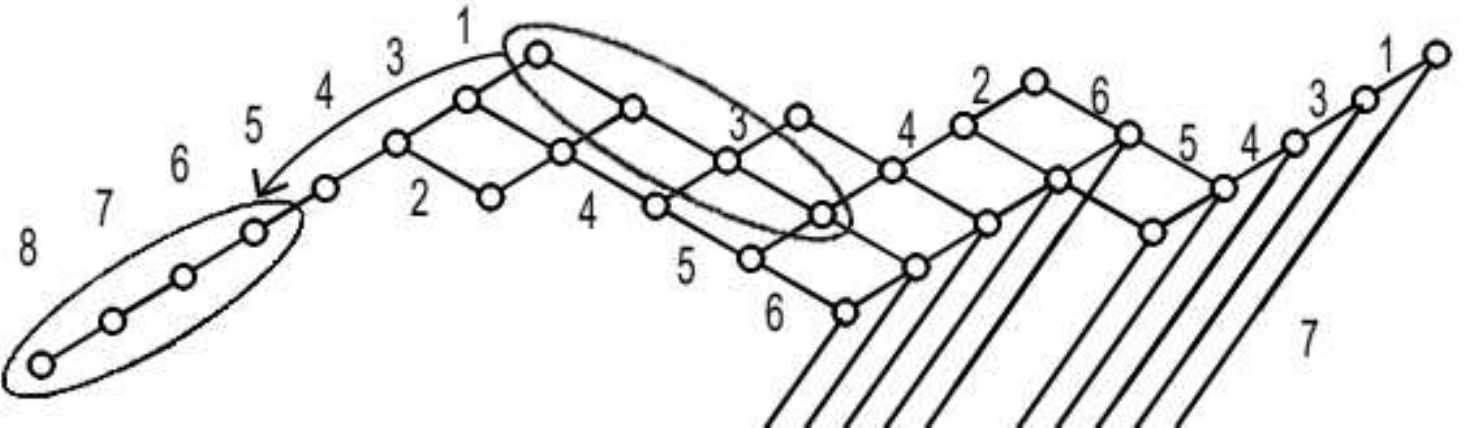}
		\end{center}
		\caption{Part of $(E_8,\vpi_8)$ and action of $w$}
		\label{PartWeightDiagramE8}
	\end{figure}

Consider the weight diagram for $(D_8,\vpi_8)$. If we cut it through the bonds marked with 2 (here marks refer to the numbering of simple root in $D_8$ as shown at Figure~\ref{NumerationDl}), then we obtain the union of diagrams $(D_6,\vpi_6)$, $(D_6,\vpi_5)\otimes (A_1,\vpi_1)$, and $(D_6,\vpi_6)$.

\begin{figure}[h]
	\begin{center}
		\includegraphics[scale=0.5]{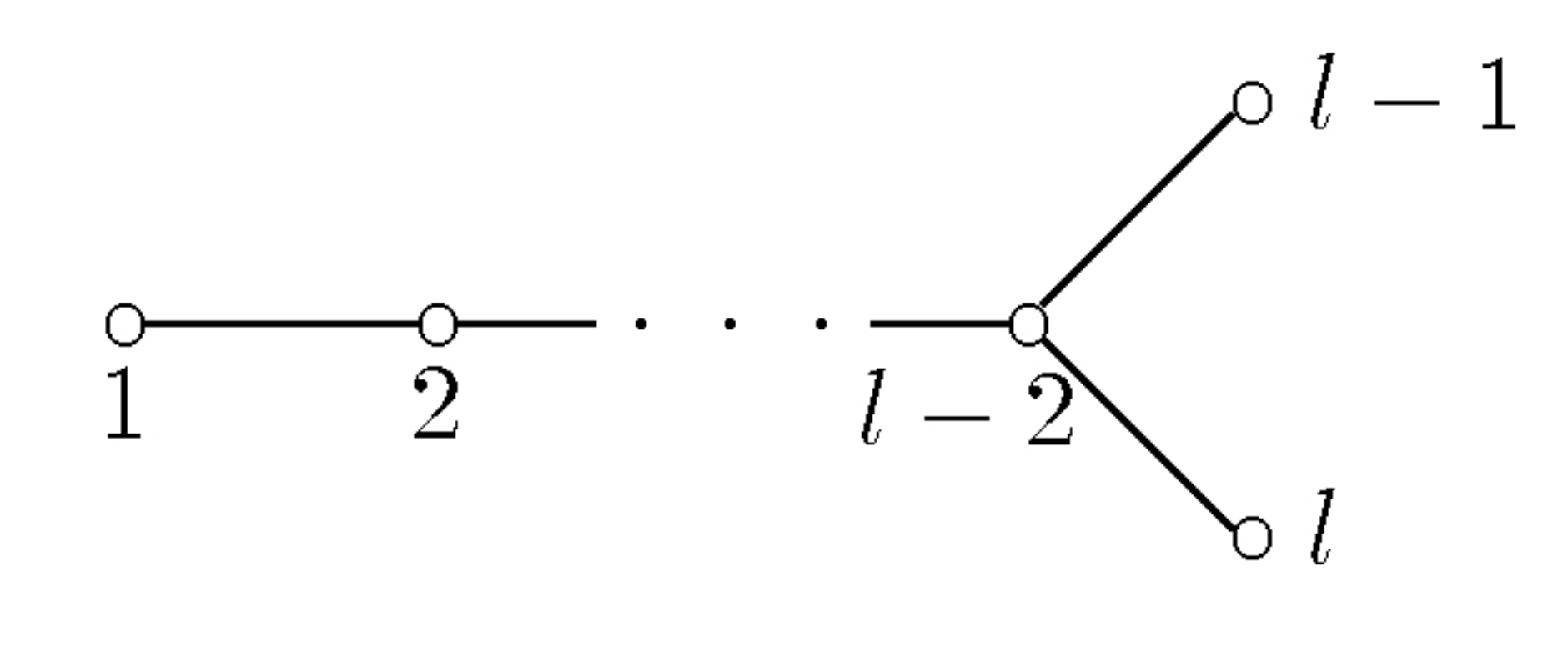}
	\end{center}
	\caption{Numbering of simple roots in $D_l$}
	\label{NumerationDl}
\end{figure}

Diagram for $(D_6,\vpi_5)$ differs from the diagram for $(D_6,\vpi_6)$ by swaping two labels; so essentially we have 4 copies of diagram $(D_6,\vpi_6)$. We give number 1 to the one containing the highest weight, number 2 to the upper half of the component $(D_6,\vpi_5)\otimes (A_1,\vpi_1)$, number 3 to its lower half, and number 4 to the one containing the lowest weight. Now we give to every vertex of the diagram $(D_8,\vpi_8)$ the number of the form $i/j$, where $i$ is the number of weight in $(D_6,\vpi_6)$ according to Figure~\ref{WeightDiagramD6}, and $j$ is the number of the copy of $(D_6,\vpi_6)$.

\begin{figure}[h]
	\begin{center}
		\includegraphics[scale=0.5]{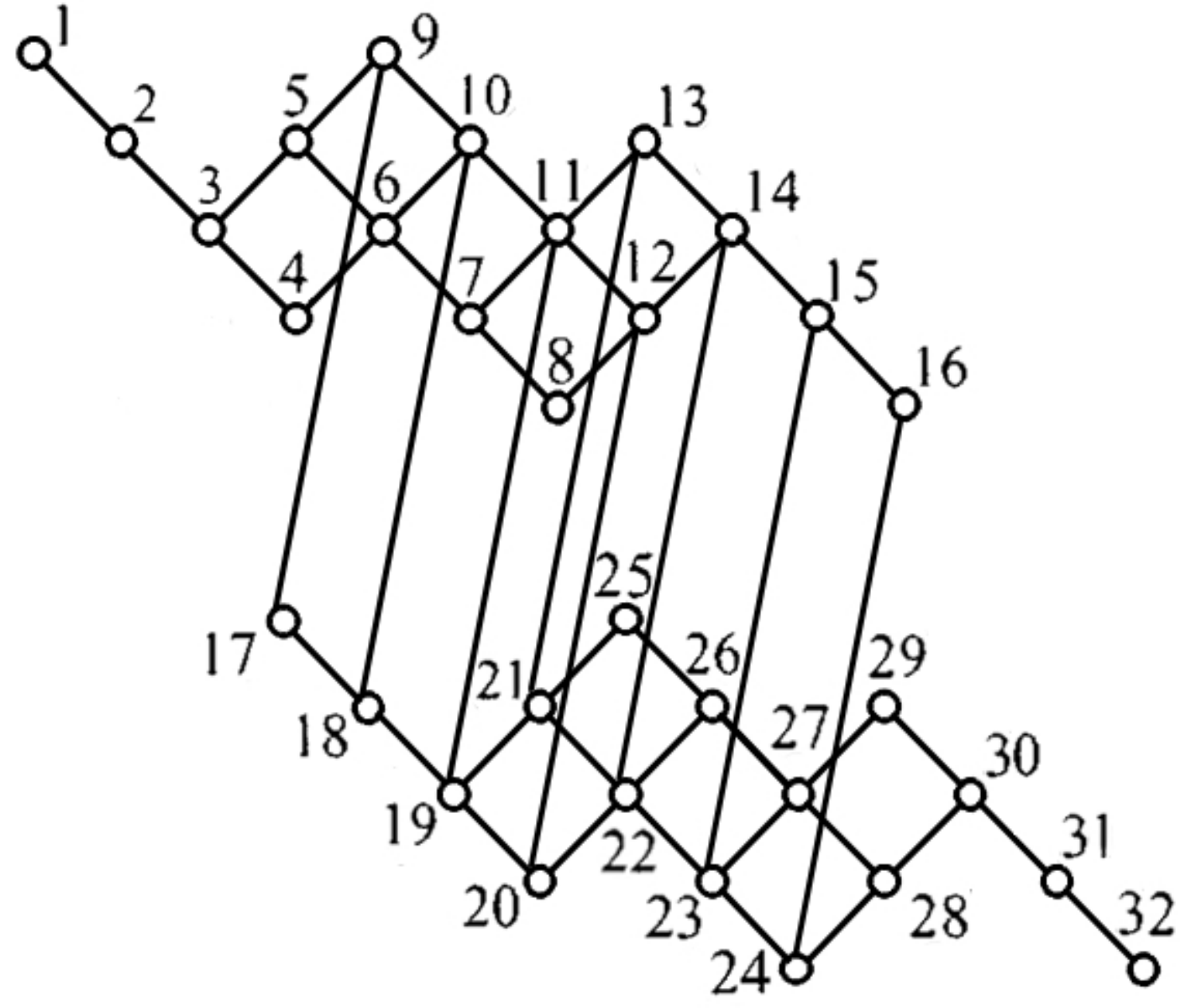}
	\end{center}
	\caption{$(D_6,\vpi_6)$}
	\label{WeightDiagramD6}
\end{figure}

Now it remains to prove the following statement. For any column vector $b=b(y)\in V_{(D_8,\vpi_8)} A$ such that it becomes unimodular in $R[y]$ and that $b_{32/4}$ is monic, there exists a column vector
$$
b^{(1)}\in E(\Phi,A)^{\le 19}b
$$
such that $b^{(1)}_{32/4}$ and $b^{(1)}_{1/1}$ are monic and the row $(b^{(1)}_{1/1},b^{(1)}_{2/1},b^{(1)}_{3/1},b^{(1)}_{5/1})$ is unimodular in $R[y]$. 

We prove this statement similarly to how we prove it for $(E_6,\vpi_1)$ and $(E_7,\vpi_7)$. Here we perform the following steps (numbering of roots is as at Figure~\ref{NumerationDl}). 

\smallskip

	{\bf Step 1.} Make the polynomial $b_{32/1}$ monic and simultaneously make the row that consists of elements $\{b_{i/j}\colon 1\le i\le 32,\, 1\le j\le 3\}\cup\{b_{32/4}\}$ unimodular in $R[y]$ by the element~$x_{\delta_{D_8}}(\xi)$.
	
\smallskip

	{\bf Step 2.} Make the polynomial $b_{8/1}$ monic and simultaneously make the row that consists of elements $\{b_{i/j}\colon 1\le i\le 24,\, 1\le j\le 3\}\cup\{b_{32/1},b_{8/4}\}$ unimodular in $R[y]$ by the element~$x_{\delta_{D_6}}(\xi)$, where $\delta_{D_6}$ is the maximal root of the subsystem generated by $\alpha_3,\ldots\alpha_8$.

\smallskip

	{\bf Step 3.} Make the row of elements $\{b_{i/j}\colon 1\le i\le 16,\, 1\le j\le 3\}\cup\{b_{32/1},b_{8/4}\}$ unimodular in $R[y]$ by the element~$x_{\alpha_3}(\xi)$ (the polynomial $b_{8/1}$ remains the same).

\smallskip

	{\bf Step 4.} Make the polynomial $b_{2/1}$ monic and simultaneously make the row that consists of elements $\{b_{i/j}\colon 1\le i\le 14,\, 1\le j\le 3\}\cup\{b_{26/1},b_{2/4}\}$ unimodular in $R[y]$ by the element~$x_{\delta_{D_4}}(\xi)$, where $\delta_{D_6}$ is the maximal root of the subsystem generated by $\alpha_5,\ldots\alpha_8$.
	
\smallskip

	{\bf Step 5.} Make the row of elements $\{b_{i/j}\colon 1\le i\le 12,\, 1\le j\le 3\}\cup\{b_{26/1},b_{2/4}\}$ unimodular in $R[y]$ by the element~$x_{\alpha_5}(\xi)$ (the polynomial $b_{2/1}$ remains the same).
	
\smallskip

	{\bf Step 6.} Make the row of elements $\{b_{i/j}\colon 1\le i\le 8,\, 1\le j\le 3\}\cup\{b_{22/1},b_{2/4}\}$ unimodular in $R[y]$ by the element~$x_{\alpha_4}(\xi)$ (the polynomial $b_{2/1}$ remains the same).
	
\smallskip
	
	{\bf Step 7.} Make the row of elements $\{b_{i/1}\colon 1\le i\le 8\}\cup\{b_{i,j}\colon 1\le i\le 7,\, 2\le j\le 3\}\cup\{b_{22/1},b_{2/4}\}$ unimodular in $R[y]$ by the element~$x_{\alpha_7}(\xi)$ (the polynomial $b_{2/1}$ remains the same).

\smallskip

	{\bf Step 8.} Make the row of elements $\{b_{i/1}\colon 1\le i\le 8\}\cup\{b_{i,j}\colon 1\le i\le 6,\, 2\le j\le 3\}\cup\{b_{22/1},b_{2/4}\}$ unimodular in $R[y]$ by the element~$x_{\alpha_6+\alpha_8}(\xi)$ (the polynomial $b_{2/1}$ remains the same).

\smallskip

{\bf Step 9.} Make the polynomial $b_{1/1}$ monic and simultaneously make the row that consists of elements $\{b_{i/1}\colon 1\le i\le 7\}\cup\{b_{i,j}\colon i\in\{1,2,3,5\},\, 2\le j\le 3\}\cup\{b_{21/1},b_{1/4}\}$ unimodular in $R[y]$ by the element~$x_{\alpha_8}(\xi)$.

\smallskip

{\bf Step 10.} Make the row $(b^{(1)}_{1/1},b^{(1)}_{2/1},b^{(1)}_{3/1},b^{(1)}_{5/1})$ unimodular in $R[y]$ by the element~$x_{\alpha_7}(\xi_{10})x_{\alpha_6}(\xi_9)x_{\alpha_4}(\xi_8)x_{\alpha_5}(\xi_7)x_{\alpha_3}(\xi_6)x_{\alpha_4}(\xi_5)x_{\alpha_2}(\xi_4)x_{\alpha_3}(\xi_3)x_{\alpha_1}(\xi_2)x_{\alpha_2}(\xi_1)$.
	\end{proof}

Now we are ready to prove Proposition~\ref{KillTheVariable}.

\smallskip

For simplicity, we will write
$$
\xymatrix{
	v\ar[r]^{N} & w
}
$$
instead of
$$
w\in E(\Phi,R)^{\le N}v\tc
$$
where $v$ and $w$ are columns in $V_{\vpi}$.

\smallskip

Set
\begin{gather*}
	N_1=\begin{cases}
		65 &\text{ for } D_5\le E_6\tc\\
		94 &\text{ for } E_6\le E_7\tc\\
		152 &\text{ for } E_7\le E_8\tp
	\end{cases}
\end{gather*} and \begin{gather*}
N_2=\begin{cases}
	7 &\text{ for } D_5\le E_6\tc\\
	10 &\text{ for } E_6\le E_7\tc\\
	139 &\text{ for } E_7\le E_8\tp
\end{cases}
\end{gather*}

Applying Lemmas~\ref{BSInduction} and~\ref{PrepareToAddToY} $d$ times, we obtain elements $s_1$,$\ldots$,$s_{d}\in B$ and columns $b=b^{(0)}$,$ b^{(1)}$,$\ldots$,$b^{(d)}\in \Um'_{\vpi} A$ such that, firstly, 
$$
\xymatrix{
	b^{(i)}\ar[rr]^{N_2} & & b^{(i+1)} & &  i=0,\ldots,d-1\tc
}
$$
secondly, $s_i\in\<b^{(i)}_{\lambda}\colon \lambda\in\Lambda_0'\>$ for $i=1,\ldots,d$, and thirdly, $\BSdim B/(s_1,\ldots s_{i+1})<\BSdim B/(s_1,\ldots,s_i)$ for $i=0,\ldots,d-1$. In particular, the elements $s_1$,$\ldots$,$s_{d}$ generate the unit ideal.

By Lemma~\ref{AddToY} we have
$$
\xymatrix{
	b^{(i)}(y)\ar[rrr]^{N_1} & & & b^{(i)}(y+s_i^{m_i}z)
}
$$
in $A[z]$.

Therefore, we have the following chain of transformations in $A[z_1,\ldots,z_{d}]$:

\begin{multline*}
	\xymatrix{
		b=b^{(0)}(y)\ar[r] & b^{(1)}(y)\ar[r] & b^{(1)}(y+s_1^{m_1}z_1)\ar[r] & b^{(2)}(y+s_1^{m_1}z_1)\ar[r] & \ldots}\\ \xymatrix{b^{(d)}(y+s_1^{m_1}z_1+\ldots+s_{d-1}^{m_{d-1}}z_{d-1})\ar[r] & b^{(d)}(y+s_1^{m_1}z_1+\ldots+s_d^{m_d}z_d)\tp
	}
\end{multline*}

Thus we have
$$
\xymatrix{
	b(y)\ar[rr]^{\!\!\!\!\!\!\!\!\!\!\!\!\!\!\!\!\!\!\!\!\!\!\!\!\!\!\!\!\!\!\!\!\!\!\!\!\!\!\!\!\!\!d(N_1+N_2)} & &  b^{(d)}(y+s_1^{m_1}z_1+\ldots+s_d^{m_d}z_d)\tp
}
$$

Since the elements $s_1$,$\ldots$,$s_{d}$ generate the unit ideal, it follows that so do the elements $s_1^{m_1}$,$\ldots$,$s_{d}^{m_{d}}$. Specializing the indeterminates $z_i$ to elements in $yB$, we make $y+s_1^{m_1}z_1+\ldots+s_d^{m_d}z_d$ equal to zero; this concludes the proof of Proposition~\ref{KillTheVariable}.







\end{document}